# HALFWAY NEW CARDINAL CHARACTERISTICS

JÖRG BRENDLE, LORENZ J. HALBEISEN, LUKAS DANIEL KLAUSNER,
MARC LISCHKA, AND SAHARON SHELAH

Abstract. Based on the well-known cardinal characteristics $\mathfrak{s}$, $\mathfrak{r}$ and $\mathfrak{i}$, we introduce nine related cardinal characteristics by using the notion of asymptotic density to characterise different intersection properties of infinite sets. We prove several bounds and consistency results, e. g. the consistency of $\mathfrak{s} < \mathfrak{s}_{1/2}$, as well as several results about possible values of $\mathfrak{i}_{1/2}$.

## 1. Introduction

This research forms part of the study of cardinal characteristics of the continuum. For a general overview of cardinal characteristics, see [Bla10], [Hal17, chapter 9] and [Vau90] as well as [BJ95]. Based on the well-known cardinal characteristics

- $\mathfrak{s} := \min\{|\mathcal{S}| \mid \mathcal{S} \subseteq [\omega]^\omega \text{ and } \forall\, X \in [\omega]^\omega\ \exists\, S \in \mathcal{S}\colon |X \cap S| = |X \setminus S| = \aleph_0\}$ (the splitting number),
- $\mathfrak{r} := \min\{|\mathcal{R}| \mid \mathcal{R} \subseteq [\omega]^\omega \text{ and } \nexists\, X \in [\omega]^\omega\ \forall\, R \in \mathcal{R}\colon |R \cap X| = |R \setminus X| = \aleph_0\}$ (the reaping number), and
- $\mathfrak{i} := \min\{|\mathcal{I}| \mid \mathcal{I} \subseteq [\omega]^\omega, \forall\, \mathcal{A} \cupdot \mathcal{B} \in \mathrm{Fin}(\mathcal{I})\colon \left|\bigcap_{A \in \mathcal{A}} A \cap \bigcap_{B \in \mathcal{B}} (\omega \setminus B)\right| = \aleph_0$ and $\mathcal{I}$ is maximal with this property$\}$ (the independence number),

we were inspired to define specialised variants of these (all of them related in some way to asymptotic density, in particular asymptotic density $1/2$) and obtained a number of bounds and consistency results for them.

We use the standard notation; in addition to $\mathfrak{s}$, $\mathfrak{r}$ and $\mathfrak{i}$ mentioned above, we will refer to a few other well-known cardinal characteristics.

Given an ideal $\mathcal{I}$ on some base set $X$, we can define four cardinal characteristics:

- the additivity number $\mathrm{add}(\mathcal{I}) := \min\{|\mathcal{A}| \mid \mathcal{A} \subseteq \mathcal{I} \text{ and } \bigcup \mathcal{A} \notin \mathcal{I}\}$,
- the covering number $\mathrm{cov}(\mathcal{I}) := \min\{|\mathcal{A}| \mid \mathcal{A} \subseteq \mathcal{I} \text{ and } \bigcup \mathcal{A} = X\}$,
- the uniformity number $\mathrm{non}(\mathcal{I}) := \min\{|Y| \mid Y \subseteq X \text{ and } Y \notin \mathcal{I}\}$, and
- the cofinality $\mathrm{cof}(\mathcal{I}) := \min\{|\mathcal{A}| \mid \mathcal{A} \subseteq \mathcal{I} \text{ and } \forall\, B \in \mathcal{I}\ \exists\, A \in \mathcal{A}\colon B \subseteq A\}$.

---

2010 *Mathematics Subject Classification.* Primary 03E17; Secondary 03E35, 03E40.
The third author was supported by the Austrian Science Fund (FWF) project P29575 "Forcing Methods: Creatures, Products and Iterations".
The fifth author was partially supported by the European Research Council grant 338821 "Model Theory and Its Applications: Dependent Classes". This is the fifth author's publication #1150.
We are grateful to the anonymous referee for suggesting numerous improvements to both the content and the presentation of this paper, and to Barnabás Farkas, Martin Goldstern, Salome Schumacher and Thilo Weinert for various helpful discussions and suggestions.

In particular, we will refer to these cardinal characteristics for

- the ideal $\mathcal{N} := \{A \subseteq 2^\omega \mid \lambda(A) = 0\}$ of *Lebesgue null sets* and
- the ideal $\mathcal{M} := \{A \subseteq \omega^\omega \mid A = \bigcup_{n\in\omega} A_n \text{ and } \forall\, n < \omega\colon A_n \text{ nowhere dense}\}$ of *meagre sets*.

Finally, we will refer to two more cardinal characteristics:

- $\mathfrak{b} := \min\{|B| \mid B \subseteq \omega^\omega \text{ and } \forall\, g \in \omega^\omega\ \exists\, f \in B\colon f \not\leq^* g\}$ (the unbounding number) and
- $\mathfrak{d} := \min\{|D| \mid D \subseteq \omega^\omega \text{ and } \forall\, g \in \omega^\omega\ \exists\, f \in D\colon g \leq^* f\}$ (the dominating number).

We will use the following concept in a few of the proofs:

**Definition 1.1.** A *chopped real* is a pair $(x, \Pi)$ where $x \in 2^\omega$ and $\Pi$ is an interval partition of $\omega$. We say a real $y \in 2^\omega$ *matches* $(x, \Pi)$ if $y{\restriction}_I = x{\restriction}_I$ for infinitely many $I \in \Pi$.

We note that the set $\mathrm{Match}(x, \Pi)$ of all reals matching $(x, \Pi)$ is a comeagre set (see [Bla10, Theorem 5.2]).

We remark that we will not rigidly distinguish between a real $r$ in $2^\omega$ and the set $R := r^{-1}(1)$, or conversely, between a subset of $\omega$ and its characteristic function.

The paper is structured as follows. In section 2, we introduce and work on several cardinal characteristics related to $\mathfrak{s}$. In section 3, we introduce and work on cardinal characteristics mostly related to $\mathfrak{r}$ and $\mathfrak{i}$, and we prove a few more results on possible values of one of them ($\mathfrak{i}_{1/2}$) in section 4. The final section 5 summarises the open questions.

## 2. Characteristics Related to $\mathfrak{s}$

Recall the following concepts from number theory.

**Definition 2.1.** For $X \in [\omega]^\omega$ and $0 < n < \omega$, define the *initial density* (of $X$ up to $n$) as

$$d_n(X) := \frac{|X \cap n|}{n}$$

and the *lower* and *upper density* of $X$ as

$$\underline{d}(X) := \liminf_{n\to\infty}(d_n(X)) \quad \text{and} \quad \bar{d}(X) := \limsup_{n\to\infty}(d_n(X)),$$

respectively. In case of convergence of $d_n(X)$, call

$$d(X) := \lim_{n\to\infty}(d_n(X))$$

the *asymptotic density* or just the *density* of $X$.

We define four relations on $[\omega]^\omega \times [\omega]^\omega$ and their associated cardinal characteristics.

**Definition 2.2.** Let $S, X \in [\omega]^\omega$. We define the following relations:

- *$S$ bisects $X$ in the limit* (or just *$S$ bisects $X$*), written as $S \mid_{1/2} X$, if
$$\lim_{n\to\infty} \frac{|S \cap X \cap n|}{|X \cap n|} = \lim_{n\to\infty} \frac{d_n(S \cap X)}{d_n(X)} = \frac{1}{2}.$$
- For $0 < \varepsilon < 1/2$, *$S$ $\varepsilon$-almost bisects $X$*, written as $S \mid_{1/2\pm\varepsilon} X$, if for all but finitely many $n < \omega$ we have
$$\frac{|S \cap X \cap n|}{|X \cap n|} = \frac{d_n(S \cap X)}{d_n(X)} \in \left(\frac{1}{2} - \varepsilon, \frac{1}{2} + \varepsilon\right).$$
- *$S$ weakly bisects $X$*, written as $S \mid^{w}_{1/2} X$, if for any $\varepsilon > 0$, for infinitely many $n < \omega$ we have
$$\frac{|S \cap X \cap n|}{|X \cap n|} = \frac{d_n(S \cap X)}{d_n(X)} \in \left(\frac{1}{2} - \varepsilon, \frac{1}{2} + \varepsilon\right).$$
- *$S$ cofinally bisects $X$*, written as $S \mid^{\infty}_{1/2} X$, if for infinitely many $n < \omega$ we have
$$\frac{|S \cap X \cap n|}{|X \cap n|} = \frac{d_n(S \cap X)}{d_n(X)} = \frac{1}{2}.$$

**Definition 2.3.** We say a family $\mathcal{S}$ of infinite sets is

$$\begin{cases} \text{bisecting (in the limit)} \\ \varepsilon\text{-almost bisecting} \\ \text{weakly bisecting} \\ \text{cofinally bisecting} \end{cases}$$

if for each $X \in [\omega]^\omega$ there is some $S \in \mathcal{S}$ such that

$$\begin{cases} S \text{ bisects } X \text{ (in the limit)} \\ S\ \varepsilon\text{-almost bisects } X \\ S \text{ weakly bisects } X \\ S \text{ cofinally bisects } X \end{cases}$$

and denote the least cardinality of such a family by $\mathfrak{s}_{1/2}$, $\mathfrak{s}_{1/2\pm\varepsilon}$, $\mathfrak{s}^{w}_{1/2}$, $\mathfrak{s}^{\infty}_{1/2}$, respectively.

**Theorem 2.4.** *The relations shown in Figure 1 hold.*

*Proof.* Recall that it is known that $\mathfrak{s} \leq \mathrm{non}(\mathcal{M})$ and $\mathfrak{s} \leq \mathrm{non}(\mathcal{N})$ (see e. g. [Bla10, Theorem 5.19]) as well as $\mathfrak{s} \leq \mathfrak{d}$ (see e. g. [Hal17, Theorem 9.4] or [Bla10, Theorem 3.3]).

$\boldsymbol{\mathfrak{s} \leq \mathfrak{s}^{w}_{1/2} \leq \mathfrak{s}^{\infty}_{1/2}}$: A cofinally bisecting real is a weakly bisecting real (being equal to $1/2$ infinitely often implies entering an arbitrary $\varepsilon$-neighbourhood of $1/2$ infinitely often), and a weakly bisecting real is a splitting real (if a real $X$ does not split another real $Y$, the relative initial density of $X$ in $Y$, that is
$$\frac{d_n(X \cap Y)}{d_n(Y)},$$

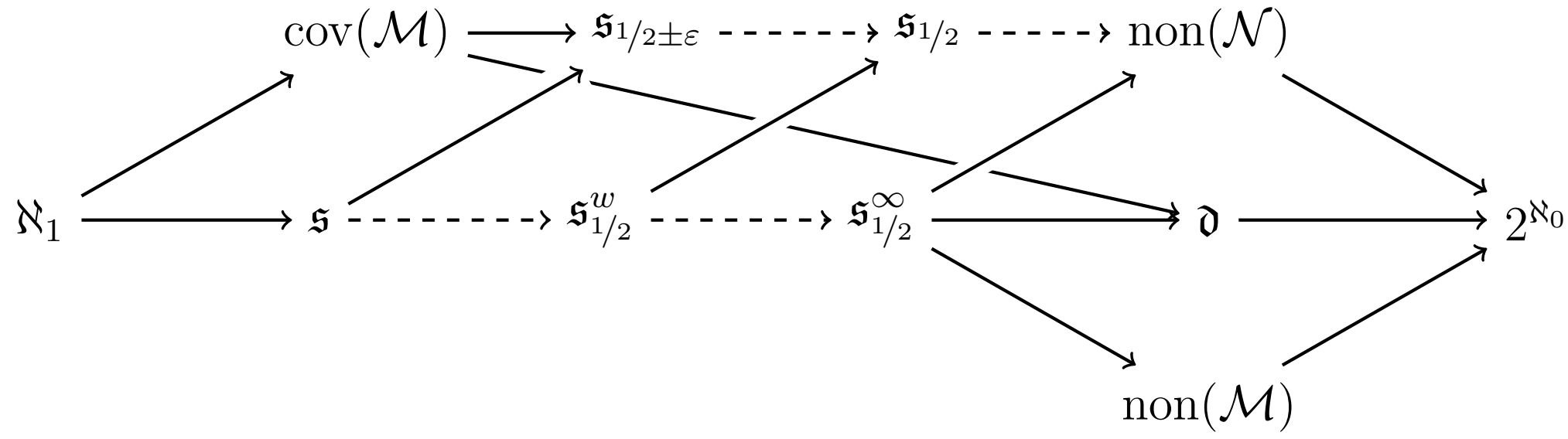

FIGURE 1. The ZFC-provable and/or consistent inequalities between $\mathfrak{s}_{1/2}$, $\mathfrak{s}_{1/2\pm\varepsilon}$, $\mathfrak{s}^w_{1/2}$, $\mathfrak{s}^\infty_{1/2}$ and other well-known cardinal characteristics, where $a \longrightarrow b$ means "$a \leq b$, consistently $a < b$" and $a \dashrightarrow b$ means "$a \leq b$, possibly $a = b$".

tends to either 0 or 1 and hence cannot be close to $1/2$ infinitely often). Hence a family witnessing the value of $\mathfrak{s}^\infty_{1/2}$ gives an upper bound for the value of $\mathfrak{s}^w_{1/2}$ (and analogously for $\mathfrak{s} \leq \mathfrak{s}^w_{1/2}$).

$\boldsymbol{\mathfrak{s} \leq \mathfrak{s}_{1/2\pm\varepsilon} \leq \mathfrak{s}_{1/2}}$: The first claim follows since an $\varepsilon$-almost bisecting real is a splitting real by the fact that finite sets have density 0 and cofinite sets have density 1, and hence if $X$ does not split $Y$, the relative initial densities of $X$ and $\omega \setminus X$ in $Y$ tend to 0 and 1, respectively (or vice versa). The second claim follows since a bisecting real is an $\varepsilon$-almost bisecting real by definition.

$\mathbf{cov}\boldsymbol{(\mathcal{M}) \leq \mathfrak{s}_{1/2\pm\varepsilon}}$: Given a family $\mathcal{S}$ witnessing the value of $\mathfrak{s}_{1/2\pm\varepsilon}$, take $S \in \mathcal{S}$.

We inductively define a chopped real $(S, \Pi)$ based on $S$ as follows: Let the first interval of the partition $\Pi$ be $I_0 = [0, \min(S)]$. Now, for any $n \in \omega$, given $m_n := \max(I_n)$, choose $m_{n+1}$ minimal such that $I_{n+1} := [m_n+1, m_{n+1}]$ contains $(n \cdot m_n)+1$ elements of $S$.

Now, any real $X$ matching this chopped real is not $\varepsilon$-almost bisected by $S$. Indeed, whenever such an $X$ is equal to $S$ on one of the intervals $I_n$, we have

$$\begin{aligned} \frac{d_{m_n}(S \cap X)}{d_{m_n}(X)} &= \frac{|S \cap X \cap (m_{n-1}+1)| + (n-1)m_{n-1} + 1}{|X \cap (m_{n-1}+1)| + (n-1)m_{n-1} + 1} \\ &\geq \frac{(n-1)m_{n-1}+1}{m_{n-1}+1+(n-1)m_{n-1}+1} = \frac{(n-1)m_{n-1}+1}{n \cdot m_{n-1}+2} \\ &\geq \frac{(n-1)m_{n-1}}{n \cdot m_{n-1}+m_{n-1}} \geq \frac{n-1}{n+1} \geq \frac{n-2}{n}. \end{aligned}$$

As $\varepsilon < 1/2$, for sufficiently large $n$ we get $1 - 2/n > 1/2 + \varepsilon$, and since such an $X$ is equal to $S$ on $I_n$ for infinitely many $n \in \omega$, $S$ does not $\varepsilon$-almost bisect $X$ in the limit.

Now, the family $M := \mathrm{Match}(S, \langle I_n \rangle_{n \in \omega})$ is a comeagre set. The family $E(S)$ of all reals which *are* $\varepsilon$-almost bisected by $S$ is a meagre set (since $E(S) \subseteq [\omega]^\omega \setminus M$ and hence its complement is a superset of a comeagre set), and $\{E(S) \mid S \in \mathcal{S}\}$ is a $2^\omega$-covering consisting of meagre sets.

$\boldsymbol{\mathfrak{s}^w_{1/2} \leq \mathfrak{s}_{1/2}}$: A bisecting real is a weakly splitting real – for the relative density to converge to $1/2$, it has to eventually be arbitrarily close to $1/2$, and hence also within an arbitrary $\varepsilon$-neighbourhood of $1/2$ infinitely often. The same argument using the families witnessing the cardinal characteristics holds.

$\mathfrak{s}^{\infty}_{1/2} \leq \mathbf{non}(\mathcal{M})$: For a given $X \in [\omega]^{\omega}$, we show that the set $B(X)$ of reals cofinally bisecting $X$ (contains and hence) is a comeagre set. For any $F \notin \mathcal{M}$, $F \cap B(X)$ is non-empty, hence it contains a real cofinally bisecting $X$.

Given $X$ as above, we define a chopped real $(R, \Pi)$ as follows: Let $f_X \colon \omega \to X$ be the ascending enumeration of $X$. For all $n \in \omega$, define intervals $J_n := [f_X(3^n), f_X(3^{n+1}))$ and let $\Pi$ consist of the intervals $I_0 := [0, f_X(1))$ and, for all $n \in \omega$, $I_{n+1} := J_{2n} \cup J_{2n+1}$. Define $R \subsetneq X$ such that for each $n \in \omega$,

$$R \cap J_{2n} = \varnothing \qquad \text{and} \qquad R \cap J_{2n+1} = X \cap J_{2n+1}\,.$$

Suppose the real $R_0$ matches $(R, \Pi)$ and is equal to $R$ on $I_{n+1}$. Let $k := \max(J_{2n}) + 1 = \min(J_{2n+1})$ and $\ell := \max(J_{2n+1}) + 1 = \min(J_{2n+2}) = \max(I_{n+1}) + 1 = \min(I_{n+2})$. Then we have:

$$\frac{d_k(R_0 \cap X)}{d_k(X)} \leq \frac{1}{3} \qquad \text{and} \qquad \frac{d_\ell(R_0 \cap X)}{d_\ell(X)} \geq \frac{2}{3}\,,$$

and

$$\frac{d_t(R_0 \cap X)}{d_t(X)}$$

is monotone increasing for $k \leq t \leq \ell$. But this implies the existence of some $t$ with $k \leq t \leq \ell$ such that

$$\frac{d_t(R_0 \cap X)}{d_t(X)} = \frac{1}{2}\,,$$

by the following argument: Let $t$ be maximal such that

$$\frac{d_t(R_0 \cap X)}{d_t(X)} \leq \frac{1}{2}\,.$$

Assume towards a contradiction that the inequality is strict; this implies that

$$\frac{|R_0 \cap X \cap t|}{|X \cap t|} \leq \frac{\lceil |X \cap t|/2 \rceil - 1}{|X \cap t|}\,,$$

but then

$$\frac{d_{t+1}(R_0 \cap X)}{d_{t+1}(X)} \leq \frac{\lceil |X \cap t|/2 \rceil}{|X \cap (t+1)|} \leq \frac{1}{2}\,,$$

contradicting the maximality of $t$. Hence $R_0$ cofinally bisects $X$.

The set of all reals matching $(R, \Pi)$ is comeagre, as required to finish the proof above.

$\mathfrak{s}^{\infty}_{1/2} \leq \mathfrak{d}$: Let $\mathcal{D}$ be a dominating family. Without loss of generality assume that every member $g$ of $\mathcal{D}$ is strictly increasing and satisfies $g(0) > 0$. Let $X \in [\omega]^{\omega}$ and let $f_X$ be its enumeration. Pick a $g_X =: g$ from $\mathcal{D}$ that dominates $f_X$ and define $G \colon \omega \to \omega$ by $G(n) := g^{(n+1)}(0)$ for every $n < \omega$. Then, for sufficiently large $n$,

$$G(n) \leq f_X(G(n)) < g(G(n)) = G(n+1).$$

Hence (for sufficiently large $n$) every interval $[G(n), G(n+1))$ contains at least one element of $X$ and at most $G(n+1) - G(n)$ many. Now iteratively define a function $\Gamma\colon \omega \to \omega$ by $\Gamma(0) := 0$, $\Gamma(1) := G(0) = g(0)$ and $\Gamma(n+1) := G\big(\sum_{k=0}^{n} \Gamma(k)\big) = G(\Sigma_n)$ and consider the interval partition with partition boundaries $\langle \Gamma(n) \mid n < \omega\rangle$; for sufficiently large $n$, every interval

$$\begin{aligned} I_n :=& \Big[\Gamma(n), \Gamma(n+1)\Big) = \Big[G\big(\sum_{k=0}^{n-1}(\Gamma(k))\big), G\big(\sum_{k=0}^{n}(\Gamma(k))\big)\Big) \\ =& \Big[G(\Sigma_{n-1}), G(\Sigma_{n-1}+1)\Big) \cup \ldots \cup \Big[G(\Sigma_{n-1}+\Gamma(n)-1), G(\Sigma_{n-1}+\Gamma(n))\Big) \end{aligned}$$

contains at least $\Gamma(n)$ many elements of $X$ and at most $\Gamma(n+1) - \Gamma(n)$ many of them.

Let $Y_X$ be the real defined as the union of every even interval, i.e. the intervals $I_{2k} = [\Gamma(2k), \Gamma(2k+1))$. We now show that $Y_X$ cofinally bisects $X$. The number of elements of $X$ which are in any interval $I_n$ is at least as large as the lower boundary of $I_n$ and $Y_X$ is defined to alternate between consecutive intervals. Consider the consecutive intervals $I_{2k}$ and $I_{2k+1}$. By definition of $Y_X$, at the endpoint of interval $I_{2k}$, we have

$$\frac{|X \cap Y_X \cap \Gamma(2k+1)|}{|X \cap \Gamma(2k+1)|} > \frac{1}{2};$$

conversely, at the endpoint of interval $I_{2k+1}$, we have (also by definition of $Y_X$, as $I_{2k+1}$ is disjoint from $Y_X$)

$$\frac{|X \cap Y_X \cap \Gamma(2k+2)|}{|X \cap \Gamma(2k+2)|} < \frac{1}{2}.$$

This means that at some point between the endpoints of $I_{2k}$ and $I_{2k+1}$, the relative initial density has to cross over from $\geq 1/2$ to $\leq 1/2$. An easy proof by contradiction shows that going from $>1/2$ to $<1/2$ in a single step is impossible. Hence the relative initial density of $Y_X$ in $X$ is $1/2$ infinitely often.

**$\mathfrak{s}^{\infty}_{1/2} \leq \mathrm{non}(\mathcal{N})$:** Given some $X \in [\omega]^\omega$ with enumerating function $f_X$ and a Lebesgue-random set $S$ (i.e. such that $\forall\, n < \omega\colon \Pr[n \in S] = 1/2$), the function $g(n) := |X \cap S \cap f_X(n)| - n/2$ defines a balanced random walk with step size $1/2$, since

$$g(n+1) - g(n) = \begin{cases} +1/2 & f_X(n) \in S, \\ -1/2 & f_X(n) \notin S. \end{cases}$$

From probability theory we know that for almost all $S$, $g(n)$ will be 0 infinitely often. Equivalently, almost surely,

$$\frac{g(n)}{n} + \frac{1}{2} = \frac{|X \cap S \cap f_X(n)|}{n}$$

will be $1/2$ infinitely often.

In other words, for any $X \in [\omega]^\omega$, the set of all $S$ not cofinally bisecting $X$ is a null set. By contraposition, for any $X \in [\omega]^\omega$, any non-null set contains a set $S$ that cofinally bisects $X$.

**$\mathfrak{s}_{1/2} \leq \mathrm{non}(\mathcal{N})$:** Let $X \in [\omega]^\omega$ and $F \notin \mathcal{N}$. Enumerating $X =: \{x_0, x_1, x_2, \dots\}$, we define functions $f_{X,n}$ and $f_X$ as follows:

$$f_{X,n} \colon [\omega]^\omega \to \{0,1\} \colon Y \longmapsto \begin{cases} 0 & x_n \notin Y \\ 1 & x_n \in Y \end{cases}$$

$$f_X \colon [\omega]^\omega \to [0,1] \colon Y \longmapsto \begin{cases} \lim_{k\to\infty} \dfrac{\sum_{n=1}^{k} f_{X,n}(Y)}{k} & \text{if the limit exists} \\ 0 & \text{otherwise} \end{cases}$$

It is clear that $\lambda(f_{X,n}^{-1}(\{1\})) = 1/2$. Hence, the $f_{X,n}$ are identically distributed random variables on the probability space $[\omega]^\omega$ with probability measure the Lebesgue measure $\lambda$. Moreover, they are independent and have finite variance. By the law of large numbers it follows that $f_X$ is almost surely equal to $1/2$, in other words $\lambda(f_X^{-1}(\{1/2\})) = 1$. This means that with

$$S_X := \{Y \in [\omega]^\omega \mid f_X(Y) = 1/2\} = \{Y \in [\omega]^\omega \mid Y \mid_{1/2} X\},$$

we have that $\lambda(S_X) = 1$ and hence $S_X \notin \mathcal{N}$. Hence $F \cap S_X \neq \varnothing$ and there is some $S \in F$ such that $S \mid_{1/2} X$. Since all this holds for any $X \in [\omega]^\omega$, we have $\mathfrak{s}_{1/2} \leq \mathrm{non}(\mathcal{N})$.

**$\mathrm{Con}(\mathrm{non}(\mathcal{M}) < \mathfrak{s}_{1/2\pm\varepsilon})$ and $\mathrm{Con}(\mathfrak{s}^\infty_{1/2} < \mathfrak{s}_{1/2\pm\varepsilon})$:** This is implied by the consistency of $\mathrm{non}(\mathcal{M}) < \mathrm{cov}(\mathcal{M})$ as witnessed by the Cohen model.

**$\mathrm{Con}(\mathfrak{s}^\infty_{1/2} < \mathrm{non}(\mathcal{M}))$, $\mathrm{Con}(\mathfrak{s}^\infty_{1/2} < \mathfrak{d})$ and $\mathrm{Con}(\mathfrak{s}^\infty_{1/2} < \mathrm{non}(\mathcal{N}))$:** In the Cohen model, we have $\aleph_1 = \mathfrak{s} = \mathfrak{s}^\infty_{1/2} = \mathrm{non}(\mathcal{M}) < \mathrm{non}(\mathcal{N}) = \mathfrak{d}$; and in the random model, we have $\aleph_1 = \mathfrak{s}^\infty_{1/2} = \mathfrak{d} < \mathrm{non}(\mathcal{M})$.

**$\mathrm{Con}(\mathrm{cov}(\mathcal{M}) < \mathfrak{s} \leq \mathfrak{s}_{1/2})$:** In the Mathias model, we have $\mathrm{cov}(\mathcal{M}) < \mathfrak{s} = 2^{\aleph_0}$, see [Hal17, Theorem 26.14] for $\mathrm{Con}(\mathrm{cov}(\mathcal{M}) < \mathfrak{h} = 2^{\aleph_0})$ and [Bla10, Theorem 6.9] for $\mathfrak{h} \leq \mathfrak{s}$. □

Finally, we remark that $\mathfrak{b}$ is incomparable with all of our newly defined cardinal characteristics. This is because in the Blass–Shelah model [BS87], $\mathfrak{s}$ is strictly above $\mathfrak{b}$ and so are all of our characteristics; and in the Laver model, $\mathrm{non}(\mathcal{N})$ is strictly below $\mathfrak{b}$ and so are all of our characteristics.

## 3. Characteristics Related to $\mathfrak{r}$ and $\mathfrak{i}$

We define a second set of properties more closely related to $\mathfrak{i}$, although $\mathfrak{s}$ does reappear in this section.

**Definition 3.1.** A set $X \in [\omega]^\omega$ is *moderate* if $\underline{d}(X) > 0$ as well as $\bar{d}(X) < 1$.[1]

**Definition 3.2.** A family $\mathcal{I}_* \subseteq [\omega]^\omega$ is *statistically independent* or *$*$-independent* if for any set $X \in \mathcal{I}_*$ we have that $X$ is moderate and for any finite subfamily

[1] Actually, it would suffice to demand $\bar{d}(X) > 0$ as well as $\underline{d}(X) < 1$, though one would have to modify a few of the subsequent proofs.

$\mathcal{E} \subseteq \mathcal{I}_*$, the following holds:

$$\lim_{n\to\infty} \left( \frac{d_n\left(\bigcap_{E\in\mathcal{E}} E\right)}{\prod_{E\in\mathcal{E}} d_n(E)} \right) = 1$$

In the case of convergence of $d_n\big(\bigcap_{E\in\mathcal{E}} E\big)$ for any finite subfamily $\mathcal{E} \subseteq \mathcal{I}_*$, this simplifies to asking for $0 < d(X) < 1$ to hold for all $X \in \mathcal{I}_*$ and

$$\prod_{E\in\mathcal{E}} d(E) = d\big(\bigcap_{E\in\mathcal{E}} E\big)$$

to hold for any finite subfamily $\mathcal{E} \subseteq \mathcal{I}_*$.

We denote the least cardinality of a maximal $*$-independent family by $\mathfrak{i}_*$.

Recall that a family $\mathcal{I}$ of subsets of $\omega$ is called *independent* if for any disjoint finite subfamilies $\mathcal{A}, \mathcal{B} \subseteq \mathcal{I}$, the set

$$\bigcap_{A\in\mathcal{A}} A \cap \bigcap_{B\in\mathcal{B}} (\omega \setminus B)$$

is infinite. Generalising this notion leads to the following definitions (which are more obviously related to the classical $\mathfrak{i}$):

**Definition 3.3.** Let $\rho \in (0,1)$. A family $\mathcal{I}_\rho \subseteq [\omega]^\omega$ is *$\rho$-independent* if for any disjoint finite subfamilies $\mathcal{A}, \mathcal{B} \subseteq \mathcal{I}_\rho$, the following holds:

$$d\left( \bigcap_{A\in\mathcal{A}} A \cap \bigcap_{B\in\mathcal{B}} (\omega \setminus B) \right) = \rho^{|\mathcal{A}|} \cdot (1-\rho)^{|\mathcal{B}|},$$

which simplifies to $= {}^1\!/\!_2{}^{|\mathcal{A}|+|\mathcal{B}|}$ in the case of $\rho = {}^1\!/\!_2$. This definition is equivalent to demanding that for any finite $\mathcal{A} \subseteq \mathcal{I}_\rho$, the following holds:

$$d\left( \bigcap_{A\in\mathcal{A}} A \right) = \rho^{|\mathcal{A}|}$$

We denote the least cardinality of a maximal $\rho$-independent family by $\mathfrak{i}_\rho$.

Recalling the definition of $\mathfrak{r}$ as the least cardinality of a family $\mathcal{R} \subseteq [\omega]^\omega$ such that no $S \in [\omega]^\omega$ splits every $R \in \mathcal{R}$, we naturally arrive at the following definition:

**Definition 3.4.** A family $\mathcal{R}_{1/2} \subseteq [\omega]^\omega$ is *${}^1\!/\!_2$-reaping* if there is no $S \in [\omega]^\omega$ bisecting all $R \in \mathcal{R}_{1/2}$. We denote the least cardinality of a ${}^1\!/\!_2$-reaping family by $\mathfrak{r}_{1/2}$.

Given the above, the natural question is: Can we define $\mathfrak{r}_*$ analogously? Consider the following definition:

**Definition 3.5.** A family $\mathcal{R}_* \subseteq [\omega]^\omega$ is *statistically reaping* or *$*$-reaping* if

$$\nexists\, S \in [\omega]^\omega \text{ moderate such that } \forall\, X \in \mathcal{R}_*\colon \lim_{n\to\infty} \left( \frac{d_n\left(S \cap X\right)}{d_n(S) \cdot d_n(X)} \right) = 1.$$

We denote the least cardinality of a $*$-reaping family by $\mathfrak{r}_*$.

The motivation for this is as follows: Considering the analogous definitions for $\mathfrak{r}$, we might call a family $\mathcal{I}$ of infinite and coinfinite sets *maximal quasi-independent* if there is no $X$ such that for all $Y \in \mathcal{I}$ we have that $X$ splits $Y$ and $X$ splits $\omega \setminus Y$ (i. e. $X$ and $Y$ are independent for all $Y \in \mathcal{I}$). It is obvious that a reaping family is also maximal quasi-independent; the converse can easily be derived by taking a maximal quasi-independent family and saturating it (without increasing its size) by adding the complements of all its sets, resulting in a reaping family. By this train of thought, it makes sense to take Definition 3.5 as the defining property of a $*$-reaping family.

Dualising the definition of $*$-reaping leads to the following, final definitions:

**Definition 3.6.** Let $S, X \in [\omega]^\omega$ with $S$ moderate. We say that $S$ *statistically splits* $X$ or $S$ *$*$-splits* $X$, written as $S \mid_* X$, if

$$\lim_{n\to\infty} \left( \frac{d_n\,(S \cap X)}{d_n(S) \cdot d_n(X)} \right) = 1.$$

**Definition 3.7.** A family $\mathcal{S}_* \subseteq [\omega]^\omega$ is *statistically splitting* or *$*$-splitting* if

$$\forall\, X \in [\omega]^\omega \;\exists\, S \in \mathcal{S}_* \text{ moderate}\colon\, S \mid_* X.$$

We denote the least cardinality of a $*$-splitting family by $\mathfrak{s}_*$.

**Theorem 3.8.** *The relations shown in Figure 2 hold.*

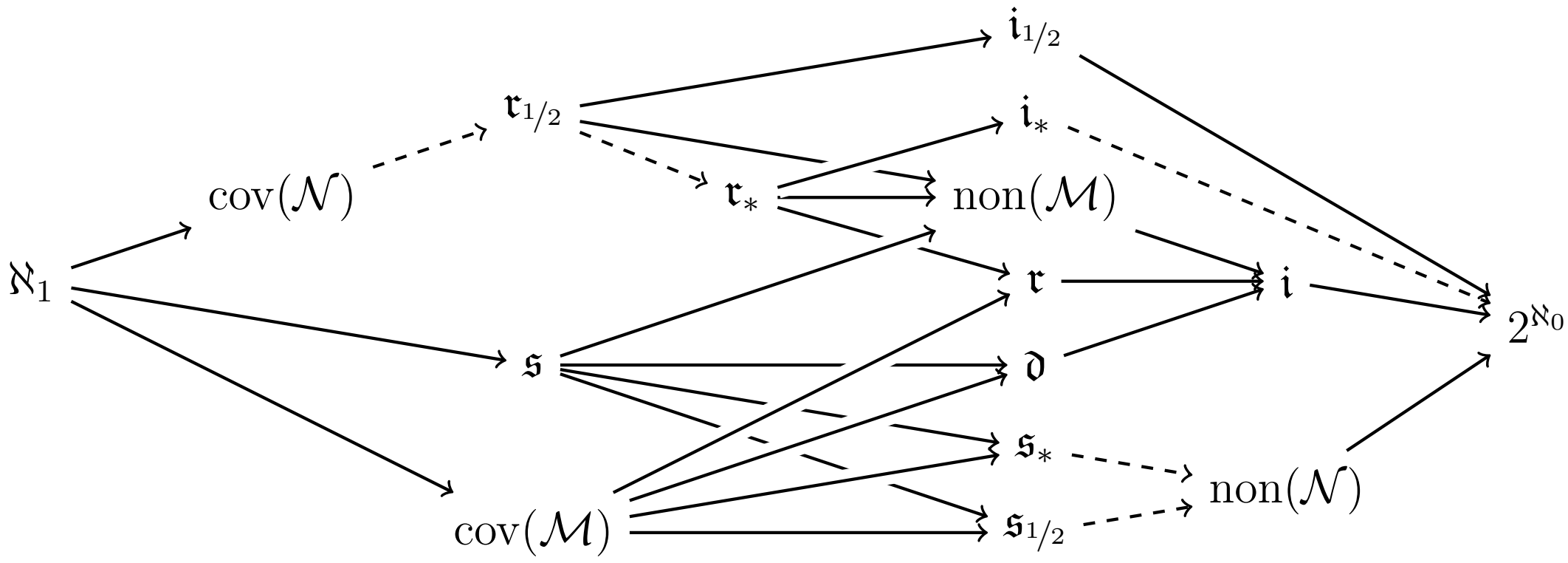

FIGURE 2. The ZFC-provable and/or consistent inequalities between $\mathfrak{i}_{1/2}$, $\mathfrak{i}_*$, $\mathfrak{r}_{1/2}$, $\mathfrak{r}_*$, $\mathfrak{s}_{1/2}$, $\mathfrak{s}_*$ and other well-known cardinal characteristics, where $a \longrightarrow b$ means "$a \leq b$, consistently $a < b$" and $a \dashrightarrow b$ means "$a \leq b$, possibly $a = b$".

*Proof.* $\mathbf{cov}(\boldsymbol{\mathcal{N}}) \leq \boldsymbol{\mathfrak{r}_{1/2}}$ and $\boldsymbol{\mathfrak{s}_*} \leq \mathbf{non}(\boldsymbol{\mathcal{N}})$: Both proofs are analogous to the proof of $\mathfrak{s}_{1/2} \leq \text{non}(\mathcal{N})$.

For the first claim, let $\mathcal{R}_{1/2}$ be a family witnessing the value of $\mathfrak{r}_{1/2}$. By the argument for $\mathfrak{s}_{1/2} \leq \text{non}(\mathcal{N})$ in the proof of Theorem 2.4, the family

$$\{[\omega]^\omega \setminus \mathcal{S}_R \mid R \in \mathcal{R}_{1/2}\}$$

is a covering of $\mathcal{N}$. (Recall that $[\omega]^\omega \setminus \mathcal{S}_R \in \mathcal{N}$ for $R \in \mathcal{R}_{1/2}$.)

For the second claim, let $X \in [\omega]^\omega$ and $F \notin \mathcal{N}$. As seen above, letting

$$S_X = \{Y \in [\omega]^\omega \mid Y \mid_{1/2} X\},$$

we have that $\lambda(S_X) = 1$ and hence $S_X \notin \mathcal{N}$. Moreover, this is true in particular for $X = \omega$ and

$$S_\omega = \{Y \in [\omega]^\omega \mid Y \mid_{1/2} \omega\} = \{Y \in [\omega]^\omega \mid d(Y) = 1/2\}.$$

Since then $F \cap S_X \cap S_\omega \neq \varnothing$, there is some $S \in F$ such that $S \mid_{1/2} X$ and $d(S) = 1/2$, which implies $S \mid_* X$.

Since all this is true for any $X \in [\omega]^\omega$, we have $\mathfrak{s}_* \leq \mathrm{non}(\mathcal{N})$.

$\boldsymbol{\mathfrak{r}_{1/2} \leq \mathfrak{r}_*}$: Let $\mathcal{R}_*$ be a $*$-reaping family and let $\mathcal{R}_{1/2} := \mathcal{R}_* \cup \{\omega\}$; clearly, $|\mathcal{R}_{1/2}| = |\mathcal{R}_*|$. Now, any $S$ which bisects all $R \in \mathcal{R}_{1/2}$ also $*$-splits all $R \in \mathcal{R}_*$ – this follows from the fact that $S \mid_{1/2} \omega$ implies $d(S) = 1/2$, and hence for any $R \in \mathcal{R}_*$, we now have

$$\frac{d_n(S \cap R)}{d_n(S) \cdot d_n(R)} = \frac{d_n(S \cap R)}{d_n(R)} \cdot \frac{1}{d_n(S)} \to 1,$$

since $S \mid_{1/2} R$ implies that the first factor converges to $1/2$, while $d(S) = 1/2$ implies that the second factor converges to 2.

$\boldsymbol{\mathfrak{r}_{1/2} \leq \mathrm{non}(\mathcal{M})}$: Since the set of all reals bisected by a fixed real $S$ is a meagre set (by the argument for $\mathrm{cov}(\mathcal{M}) \leq \mathfrak{s}_{1/2\pm\varepsilon}$), a non-meagre set contains some real not bisected by $S$ and hence is $1/2$-reaping.

$\boldsymbol{\mathfrak{r}_* \leq \mathrm{non}(\mathcal{M})}$: This is analogous to the proof of $\mathfrak{r}_{1/2} \leq \mathrm{non}(\mathcal{M})$, since the set of all reals $*$-split by a fixed moderate real $S$ is a meagre set, as well. To see this, iteratively define a chopped real based on $S$ as follows: Letting $f_S$ be the increasing enumeration of $S$, define an interval partition by defining $I_0 := [0, f_S(1)]$ and defining $I_{n+1}$ such that it contains $2^n \cdot \max(I_n)$ elements of $S$. The sets matching this chopped real form a comeagre set which consists of reals $X$ not $*$-split by $S$: As the matching intervals grow longer and longer, they "pull" $\frac{d_n(S \cap X)}{d_n(X)}$ above $1 - 1/2^n$,[2] which implies that $\frac{d_n(S \cap X)}{d_n(S) \cdot d_n(X)}$ cannot converge to 1 as $d_n(S)$ does not converge to 1 by the moderacy of $S$.

$\boldsymbol{\mathrm{cov}(\mathcal{M}) \leq \mathfrak{s}_*}$: This is analogous to the proof of $\mathrm{cov}(\mathcal{M}) \leq \mathfrak{s}_{1/2\pm\varepsilon}$ by the same argument as in the proof of $\mathfrak{r}_* \leq \mathrm{non}(\mathcal{M})$.

$\boldsymbol{\mathfrak{s} \leq \mathfrak{s}_*}$: Let $\mathcal{S}_*$ be a family witnessing the value of $\mathfrak{s}_*$ and let $X \in [\omega]^\omega$ be arbitrary. We will prove by contradiction that there must be some $S \in \mathcal{S}_*$ splitting $X$. Suppose not, that is, suppose that for any $S \in \mathcal{S}_*$, either (a) $S \cap X$ is finite or (b) $X \setminus S$ is finite. In case (a), we use the fact that $S$ is moderate to see that $d_n(S)$ must eventually be bounded from below by some $\varepsilon$, and the fact that $S \cap X$ is finite to see that $|S \cap X \cap n|$ is bounded by some $k^*$. Letting $k_n := |X \cap n|$, this eventually yields

$$\frac{d_n\,(S \cap X)}{d_n(S) \cdot d_n(X)} \leq \frac{k^*/n}{\varepsilon \cdot k_n/n} = \frac{k^*}{\varepsilon \cdot k_n} \to 0.$$

Similarly, in case (b) we use the moderacy of $S$ to see that $d_n(S)$ is eventually bounded from above by some $1 - \delta$, and the fact that $X \setminus S$ is finite to see that

[2] The strict argument for this claim is analogous to the proof of $\mathrm{cov}(\mathcal{M}) \leq \mathfrak{s}_{1/2\pm\varepsilon}$ in Theorem 2.4.

$|S \cap X \cap n|$ is bounded from below by $k_n - k^*$ for some $k^*$. (This bound simply states that after some finite aberrations, $S$ contains all elements of $X$.) Taken together, we eventually have

$$\frac{d_n\,(S \cap X)}{d_n(S) \cdot d_n(X)} \geq \frac{(k_n - k^*)/n}{(1-\delta) \cdot k_n/n}$$
$$= \frac{1}{1-\delta} - \frac{k^*}{(1-\delta) \cdot k_n} \to \frac{1}{1-\delta} = 1 + \varepsilon$$

for some $\varepsilon > 0$. In summary, for all $S \in \mathcal{S}_*$ we have that $S$ does not $*$-split $X$, and hence $\mathcal{S}_*$ could not have been a witness for the value of $\mathfrak{s}_*$.

**$\mathfrak{r}_* \leq \mathfrak{r}$**: In the previous proof, we have already shown that $*$-splitting implies splitting; this in turn proves $\mathfrak{r}_* \leq \mathfrak{r}$.

**$\mathfrak{r}_{1/2} \leq \mathfrak{i}_{1/2}$ and $\mathfrak{r}_* \leq \mathfrak{i}_*$**: For the first claim, let $\mathcal{I}_{1/2}$ be a maximal $1/2$-independent family. Define

$$\mathcal{R}_{1/2} := \left\{ \bigcap_{A \in \mathcal{A}} A \;\middle|\; \mathcal{A} \subseteq \mathcal{I}_{1/2} \right\}.$$

Then $\mathcal{R}_{1/2}$ is a $1/2$-reaping family, since the existence of an $S \in [\omega]^\omega$ bisecting each $R \in \mathcal{R}_{1/2}$ (in the limit) would contradict the maximality of $\mathcal{I}_{1/2}$.

The proof of the second claim is analogous: Take all finite tuples of sets in the witness $\mathcal{I}_*$ of the value of $\mathfrak{i}_*$ and collect their intersections in a family $\mathcal{R}_*$; this family must then be $*$-reaping, because a set $S$ $*$-splitting each $R \in \mathcal{R}_*$ would violate the maximality of $\mathcal{I}_*$, and thus $\mathcal{R}_*$ witnesses $\mathfrak{r}_* \leq \mathfrak{i}_*$.

**$\mathfrak{i}_\rho \leq 2^{\aleph_0}$ and $\mathfrak{i}_* \leq 2^{\aleph_0}$**: For $\mathfrak{i}_\rho$, consider the collection $\mathcal{I}_\rho$ of all $\rho$-independent families. Now, $\mathcal{I}_\rho$ has finite character, i. e. for each $I \subseteq 2^{\aleph_0}$, $I$ belongs to $\mathcal{I}_\rho$ iff every finite subset of $I$ belongs to $\mathcal{I}_\rho$. Hence we can apply Tukey's lemma and see that $\mathcal{I}_\rho$ has a maximal element with respect to inclusion. Therefore, $\mathfrak{i}_\rho$ is well defined and hence $\mathfrak{i}_\rho \leq 2^{\aleph_0}$. The proof for $\mathfrak{i}_*$ is analogous.

**Con($\mathfrak{r}_* < \mathfrak{r}$)**: This follows from Con(non($\mathcal{M}$) $<$ cov($\mathcal{M}$)), but we also have an explicit proof of this.

We will show that Cohen forcing does not increase $\mathfrak{r}_*$ due to the ground model reals remaining $*$-reaping; we already know that Cohen forcing increases $\mathfrak{r}$, proving our consistency statement.

Let $\dot{X}$ be a $\mathbb{C}$-name for a moderate real. We will construct a ground model real $Y$ such that for any $q \in \mathbb{C}$, we can find $r \leq q$ such that $r \Vdash \dot{X} \not\mid_* Y$.

Now pick an enumeration $\langle p_k \mid k < \omega \rangle$ of $\mathbb{C}$ which enumerates each element infinitely often. In the following argument, for each $k < \omega$, let $L_k := \sum_{\ell \leq k} \ell_k$.

- For $k = 0$, we find $q_0 \leq p_0$, $\ell_0 \geq 2$ and $A_0 \subseteq [0, \ell_0)$ such that $q_0$ decides $\dot{X}{\restriction}_{\ell_0}$, $q_0 \Vdash \dot{X}{\restriction}_{\ell_0} = A_0$ and such that $|A_0| \geq 1$, $|[0, \ell_0) \setminus A_0| \geq 1$, and at least one of these two inequalities is an equality.
- For $0 < k < \omega$, we find $q_k \leq p_k$, $\ell_k < \omega$ and $A_k \subseteq [L_{k-1}, L_k)$ such that $q_k$ decides $\dot{X}{\restriction}_{L_k}$, $q_k \Vdash \dot{X}{\restriction}_{[L_{k-1}, L_k)} = A_k$ and such that $|A_k| \geq 3L_{k-1}$, $|[L_{k-1}, L_k) \setminus A_k| \geq 3L_{k-1}$, and at least one of these inequalities is an equality.

Define $Y$ piecewise by $Y\restriction_{[L_{k-1},L_k)} := [L_{k-1}, L_k) \setminus A_k$. Assume $\dot{X}$ $*$-splits $Y$; then there must be some $q \in \mathbb{C}$ forcing this. In particular, this means that $q$ forces that for any $\varepsilon > 0$, there is some $m_\varepsilon < \omega$ such that for any $j > m_\varepsilon$,

$$\frac{d_j(\dot{X} \cap Y)}{d_j(\dot{X}) \cdot d_j(Y)} > 1 - \varepsilon.$$

Pick some sufficiently small $\varepsilon$, say $\varepsilon := 2/9$, strengthen $q$ to decide the value of $m_{2/9}$, and find $n < \omega$ such that $p_n = q$ and $L_n > m_{2/9}$. Letting $I_n := |A_n|$ and $O_n := \ell_n - I_n$, $q_n \leq q$ forces

$$d_{L_n}(\dot{X} \cap Y) \leq \frac{L_{n-1}}{L_n},$$
$$d_{L_n}(\dot{X}) \geq \frac{I_n}{L_n},$$
$$d_{L_n}(Y) \geq \frac{O_n}{L_n}.$$

Without loss of generality, $O_n = 3L_{n-1}$ and $I_n = 3L_{n-1} + \Delta$ for some $\Delta < \omega$. Then $q_n$ forces

$$\frac{d_{L_n}(\dot{X} \cap Y)}{d_{L_n}(\dot{X}) \cdot d_{L_n}(\dot{Y})} \leq \frac{\frac{L_{n-1}}{L_n}}{\frac{O_n I_n}{L_n^2}} = \frac{L_{n-1}L_n}{O_n I_n} = \frac{L_{n-1}(L_{n-1} + O_n + I_n)}{O_n I_n}$$
$$= \frac{L_{n-1}(7L_{n-1} + \Delta)}{3L_{n-1}(3L_{n-1} + \Delta)} = \frac{7L_{n-1} + \Delta}{3 \cdot (3L_{n-1} + \Delta)},$$

which is strictly decreasing in $\Delta$ and is $7/9$ for $\Delta = 0$. This contradicts the assumption on $q$, proving that $\dot{X}$ does not $*$-split $Y$ in $V^{\mathbb{C}}$.

Hence assuming CH in the ground model and forcing with $\mathbb{C}_\lambda$ for some $\lambda \geq \aleph_2$ with $\lambda = \lambda^{\aleph_0}$ gives us $V^{\mathbb{C}_\lambda} \vDash \mathfrak{r}_* = \aleph_1 < \lambda = \mathfrak{r} = \mathfrak{c}$.

**Con($\mathfrak{r}_{1/2} <$ non($\mathcal{M}$)) and Con($\mathfrak{r}_* <$ non($\mathcal{M}$))**: This follows from Con($\mathfrak{r} <$ non($\mathcal{M}$)), see [BJ95, Model 7.5.9].

**Con($\mathfrak{s} < \mathfrak{s}_*$)**: Just like Con($\mathfrak{r}_* < \mathfrak{r}$), this follows from Con(non($\mathcal{M}$) $<$ cov($\mathcal{M}$)), but once more, we also have an explicit proof of this.

We will show that Cohen forcing increases $\mathfrak{s}_*$ due to the Cohen real not being $*$-split by any real from the ground model; we already know that Cohen forcing keeps $\mathfrak{s}$ small, proving our consistency statement.

The proof uses the same technique as the one for $\mathfrak{s} \leq \mathfrak{s}_*$: Given some moderate $X \in [\omega]^\omega \cap V$, with moderacy in the sense of $\bar{d}(X) = 1 - 2\varepsilon$ and $d_n(X) < 1 - \varepsilon$ for all $n \geq n_0$ for some $n_0$, we will show that the assumption that there is a condition forcing $X \mid_* \dot{C}$, i.e. that $X$ $*$-splits the Cohen real, leads to a contradiction.

So suppose that there were some $p \in \mathbb{C}$ such that $p \Vdash X \mid_* \dot{C}$; more specifically, suppose that for some $n_1$, even $p \Vdash \frac{d_n(X \cap \dot{C})}{d_n(X) \cdot d_n(\dot{C})} < 1 + \delta$ for all $n \geq n_1$, where $\delta := \frac{\varepsilon/2}{1-\varepsilon}$.

We now define $q \leq p$ as follows: Let $n_2$ be large enough such that

$$\frac{|p|}{|X \cap [|p|, n_2)|} < \frac{\varepsilon}{2} \quad \Longleftrightarrow \quad \frac{2 \cdot |p|}{\varepsilon} < |X \cap [|p|, n_2)|;$$

this is possible due to the moderacy of $X$ (which implies $X$ is infinite). Let $k := \max(n_0, n_1, n_2)$ and $q := p^\frown \chi_X{\restriction}_{[|p|,k)}$, that is, extend $p$ by the next $k - |p|$ values of the characteristic function of $X$. Then we have

$$\frac{d_k(X \cap \dot{C})}{d_k(X) \cdot d_k(\dot{C})} > \frac{1}{1-\varepsilon} \cdot \frac{d_k(X \cap \dot{C})}{d_k(\dot{C})}$$

by the moderacy of $X$ and by $k \geq n_0$. By our choice of $q$, we have

$$\begin{aligned} q \Vdash \frac{d_k(X \cap \dot{C})}{d_k(\dot{C})} &= \frac{|X \cap \dot{C} \cap k|}{|\dot{C} \cap k|} \geq \frac{|\dot{C} \cap k| - |p|}{|\dot{C} \cap k|} = 1 - \frac{|p|}{|\dot{C} \cap k|} \\ &\geq 1 - \frac{|p|}{|\dot{C} \cap [|p|, k)|} = 1 - \frac{|p|}{|X \cap [|p|, k)|} > 1 - \frac{\varepsilon}{2}, \end{aligned}$$

with the first inequality being an equality in the "worst case" of $X{\restriction}_{|p|} \equiv 0$ and $(p = q{\restriction}_{|p|} =)\ \dot{C}{\restriction}_{|p|} \equiv 1$ and the final inequality following from $k \geq n_2$. This implies that

$$q \Vdash \frac{d_k(X \cap \dot{C})}{d_k(X) \cdot d_k(\dot{C})} > \frac{1 - \varepsilon/2}{1 - \varepsilon} = 1 + \delta,$$

contradicting (since $k \geq n_1$) the original assumption on $p$.

**Con(cov($\mathcal{M}$) $< \mathfrak{s} \leq \mathfrak{s}_*$):** Follows as in the proof of Con(cov($\mathcal{M}$) $< \mathfrak{s} \leq \mathfrak{s}_{1/2}$).

**Con($\mathfrak{r}_{1/2} < \mathfrak{i}_{1/2}$)** and **Con($\mathfrak{r}_* < \mathfrak{i}_*$):** See Lemma 3.9 and Corollary 3.13 below.

**Con($\mathfrak{i}_{1/2} < 2^{\aleph_0}$):** This follows from Lemma 3.16 below. ☐

**Lemma 3.9.** Con($\mathfrak{r}_{1/2} < \mathfrak{i}_{1/2}$).

We will prove the following: Assume CH in the ground model and let $\lambda > \mu > \aleph_1$ be regular cardinals with $\lambda = \lambda^{\aleph_0}$. Then there is a forcing extension satisfying $\mathrm{add}(\mathcal{N}) = \mathrm{cof}(\mathcal{N}) = \mathfrak{r}_{1/2} = \mu$ and $\mathfrak{c} = \mathfrak{i}_{1/2} = \lambda$.

For the proof, we shall use the method of templates, originally due to the fifth author [She04]. More explicitly, our model is the model from [Bre02, Proposition 4.7], that is, the original template model with localisation forcing instead of Hechler forcing. While we try to be as self-contained as possible and – in particular – provide the combinatorial argument in full detail below, we cannot repeat all basic lemmata of the technique and therefore assume some knowledge of the method. We recommend that the reader unfamiliar with templates works with either the survey paper [Bre02] or, preferably, the more recent [Bre21].

Given a linear order $\langle L, \leq_L \rangle$ and $x \in L$, let $L_x = \{y \in L \mid y <_L x\}$ be the initial segment determined by $x$.

**Definition 3.10** ([Bre21, Definition 21], see also [Bre02, p. 8])**.** A *template* is a pair $(L, \bar{\mathcal{I}})$ with $\bar{\mathcal{I}} = \{\mathcal{I}_x \mid x \in L\}$ such that $\langle L, \leq_L \rangle$ is a linear order, $\mathcal{I}_x \subseteq \mathcal{P}(L_x)$ for all $x \in L$, and

(1) $\mathcal{I}_x$ contains all singletons and is closed under unions and intersections,
(2) $\mathcal{I}_x \subseteq \mathcal{I}_y$ for $x <_L y$, and
(3) $\mathcal{I} := \bigcup_{x \in L} \mathcal{I}_x \cup \{L\}$ is well-founded with respect to inclusion, as witnessed by the *depth function* $\mathrm{dp}_{\mathcal{I}} \colon \mathcal{I} \to \mathrm{Ord}$.

If $A \subseteq L$ and $x \in L$, we define $\mathcal{I}_x{\restriction}_A := \{B \cap A \mid B \in \mathcal{I}_x\}$, the *trace* of $\mathcal{I}_x$ on $A$, and let $\bar{\mathcal{I}}{\restriction}_A = \{\mathcal{I}_x{\restriction}_A \mid x \in A\}$. Note that $(A, \bar{\mathcal{I}}{\restriction}_A)$ is a template as well.

Next recall (see [Bre02, p. 27]) that *localisation forcing* $\mathbb{LOC}$ consists of pairs $(\sigma, \varphi) \in ([\omega]^{<\omega})^{<\omega} \times ([\omega]^{<\omega})^{\omega}$ such that, letting $n := |\sigma|$, we have that $|\sigma(i)| \leq i$ for all $i < n$ and $|\varphi(i)| \leq n$ for all $i$. The order is given by $(\tau, \psi) \leq (\sigma, \varphi)$ if $n := |\tau| \geq m := |\sigma|$, $\tau \supseteq \sigma$, $\varphi(j) \subseteq \tau(j)$ for all $m \leq j < n$ and $\varphi(j) \subseteq \psi(j)$ for all $j$. $\mathbb{LOC}$ generically adds a *slalom*, that is, a function $\phi\colon \omega \to [\omega]^{<\omega}$ with $|\phi(n)| \leq n$ for all $n$, which *localises* all ground model reals, which means that for all $f \in \omega^\omega$ from the ground model and almost all $n$, $f(n) \in \phi(n)$. By Bartoszyński's characterisation of $\mathrm{add}(\mathcal{N})$ and $\mathrm{cof}(\mathcal{N})$ [BJ95], any iteration adding a cofinal sequence of length $\mu$ of $\mathbb{LOC}$-generic reals, where $\mu$ is uncountable regular, will force $\mathrm{add}(\mathcal{N}) = \mathrm{cof}(\mathcal{N}) = \mu$.

Since $\mathbb{LOC}$ is a correctness-preserving $\sigma$-linked forcing notion, it can be iterated along a template, see [Bre05] (see also [Bre21, Definition and Theorem 23]). This means that given a template $(L, \bar{\mathcal{I}})$ and a set $A \subseteq L$, we can define the partial order $\mathbb{P}{\restriction}_A$ by producing $\mathbb{P}{\restriction}_B$ for $B \in \mathcal{I}{\restriction}_A$ by recursion on $\mathrm{dp}_{\mathcal{I}{\restriction}_A}$ in such a way that $\mathbb{P}{\restriction}_B$ completely embeds into $\mathbb{P}{\restriction}_A$ for any $B \subseteq A \subseteq L$. "Successor" steps of the iteration are two-step iterations of the form $\mathbb{P}{\restriction}_A \cong \mathbb{P}{\restriction}_B \star \dot{\mathbb{LOC}}$ while in "limit" steps, we take direct limits. See, again, the statement and the proof of [Bre21, Definition and Theorem 23] for details. The whole iteration, $\mathbb{P}{\restriction}_L$, is still ccc by [Bre21, Lemma 24]. The proof that $\mathbb{LOC}$ is correctness-preserving – and thus fits into this framework – is exactly like the corresponding proof for Hechler forcing [Bre21, Lemma 28].[3]

We finally recall the definition of the concrete template $(L, \bar{\mathcal{I}})$ we are using here (see [Bre21, between Lemmata 28 and 29] or [Bre02, p. 20]). Let $\mu$ and $\lambda$ be cardinals. As usual, $\lambda^*$ denotes (a disjoint copy of) $\lambda$ with the reverse ordering. Elements of $\lambda$ will be called *positive* and elements of $\lambda^*$ *negative*. Choose a partition $\lambda^* = \bigcup_{\alpha<\omega_1} S^\alpha$ such that each $S^\alpha$ is coinitial in $\lambda^*$.

**Definition 3.11.** The linear order $L$ consists of non-empty finite sequences $x$ such that $x(0) \in \mu$ and $x(n) \in \lambda^* \cup \lambda$ for $n > 0$. The order is naturally given by $x < y$ if

- either $x \subsetneq y$ and $y(|x|) \in \lambda$,
- or $y \subsetneq x$ and $x(|y|) \in \lambda^*$,
- or $x(0) < y(0)$,
- or, letting $n := \min\{m \mid x(m) \neq y(m)\} > 0$, $x(n) < y(n)$ in the natural ordering of $\lambda^* \cup \lambda$.

It is immediate that this is indeed a linear order. We identify sequences of length one with their range so that $\mu \subseteq L$ is cofinal.

We write $\gamma^*$ for the element of $\lambda^*$ corresponding to $\gamma \in \lambda$. We also define $\mathrm{abs}(z) \in \lambda$ for any $z \in \lambda \cup \lambda^*$ by $\mathrm{abs}(z) = z$ for $z \in \lambda$ and $z = \mathrm{abs}(z)^*$ for $z \in \lambda^*$.

We call $x \in L$ *relevant* if $|x| \geq 3$ is odd, $x(n)$ is negative for odd $n$ and positive for even $n$, $x(|x|-1) < \omega_1$, and whenever $n < m$ are even such that $x(n), x(m) < \omega_1$, then there are $\beta < \alpha$ such that $x(n-1) \in S^\alpha$ and $x(m-1) \in S^\beta$. For relevant $x$, set $J_x = [x{\restriction}_{|x|-1}, x)$, the interval of nodes between $x{\restriction}_{|x|-1}$ and $x$ in the ordering of

[3] For an alternative explanation as to why $\mathbb{LOC}$ can be iterated along a template, see [Bre02, Lemma 4.4].

$L$. Notice that if $x < y$ are relevant, then either $J_x \cap J_y = \varnothing$ or $J_x \subsetneq J_y$ (in which case we also have $|y| \le |x|$, $x\restriction_{|y|-1} = y\restriction_{|y|-1}$ and $x(|y|-1) \le y(|y|-1)$).

**Definition 3.12.** For $x \in L$, let $\mathcal{I}_x$ consist of finite unions (1) of $L_\alpha$, where $\alpha \le x$ and $\alpha \in \mu$, (2) of $J_y$, where $y \le x$ is relevant, and (3) of $\varnothing$ and singletons.

Then $(L, \bar{\mathcal{I}})$ (with $\mathcal{I} = \{\mathcal{I}_x \mid x \in L\}$)) is indeed a template ([Bre21, Lemma 29], see also [Bre02, Lemma 3.2]). Note that, ordered by inclusion, $L$ is a tree of countable height. Countable subtrees $A, B \subseteq L$ are called *isomorphic* if there is a bijection $\varphi = \varphi_{A,B} \colon A \to B$ such that for all $x, y \in A$ and all $n \in \omega$,

- $|\varphi(x)| = |x|$,
- $\varphi(x)\restriction_n = \varphi(x\restriction_n)$,
- $x < y$ iff $\varphi(x) < \varphi(y)$,
- $x(n)$ is positive iff $\varphi(x)(n)$ is positive, and
- $\varphi$ maps $\mathcal{I}\restriction_A$ to $\mathcal{I}\restriction_B$.

Since the trace of $\mathcal{I}$ on any countable set is countable, there are at most $\mathfrak{c}$ many isomorphism types of trees. Note that, in view of the last two clauses, if $A$ and $B$ are isomorphic, then so are $\mathbb{P}\restriction_A$ and $\mathbb{P}\restriction_B$, since the partial order only depends on the structure of the template. If only the first four clauses hold, we call the trees *weakly isomorphic.*

*Proof of Lemma 3.9.* We closely follow [Bre21, Theorem 30], but provide all the details.

Take the template $(L, \bar{\mathcal{I}})$ introduced above. Let $\mathbb{P} = \mathbb{P}\restriction_L$ be the iteration of localisation forcing along this template, that is, for each $x \in L$ we have a copy of $\mathbb{LOC}$ as an iterand (see, again, [Bre21, Definition and Theorem 23] for details). It is well-known (see the proof of [Bre21, Theorem 30]) that $|\mathbb{P}| = \lambda^\omega = \lambda$ and that $\mathbb{P}$ thus forces $\mathfrak{c} \le \lambda$.

Also, letting $\dot{\phi}_\alpha$ (for $\alpha < \mu$) be the $\mathbb{P}$-name for the localisation generic added at stage $\alpha$, we see that the $\dot{\phi}_\alpha$ form a cofinal sequence of length $\mu$ of $\mathbb{LOC}$-generic reals. This is the case because $L_\alpha \in \mathcal{I}_\alpha$ and therefore $\dot{\phi}_\alpha$ is $\mathbb{LOC}$-generic over the $\mathbb{P}\restriction_{L_\alpha}$-extension. As remarked above, this implies $\mathrm{add}(\mathcal{N}) = \mathrm{cof}(\mathcal{N}) = \mu$.[4] Since we know that $\mathrm{add}(\mathcal{N}) \le \mathrm{cov}(\mathcal{N}) \le \mathfrak{r}_{1/2} \le \mathrm{non}(\mathcal{M}) \le \mathrm{cof}(\mathcal{N})$, we also have $\mathfrak{r}_{1/2} = \mu$.

We thus only have to show that $\mathfrak{i}_{1/2} \ge \lambda$. Since $\mathfrak{r}_{1/2} \le \mathfrak{i}_{1/2}$ in ZFC, we already know $\mathfrak{i}_{1/2} \ge \mu$. Thus let $\dot{\mathcal{A}}$ be a name for a $1/2$-independent family of size $<\lambda$ and $\ge\mu$, say $\dot{\mathcal{A}} = \{\dot{A}^\alpha \mid \alpha < \kappa\}$ where $\kappa \ge \omega_2 \cdot 2$ is an ordinal (chosen this way for later pruning arguments). By [Bre21, Lemma 25], there are countable $B^\alpha \subseteq L$ such that the $\dot{A}^\alpha$ are $\mathbb{P}\restriction_{B^\alpha}$-names. More explicitly, letting $\{p^\alpha_{n,i} \mid n \in \omega\}$ (for $i \in \omega$) be maximal antichains and $\{k^\alpha_{n,i} \in \{0,1\} \mid i, n \in \omega\}$ be such that $p^\alpha_{n,i} \Vdash i \in \dot{A}^\alpha$ iff $k^\alpha_{n,i} = 1$ and $p^\alpha_{n,i} \Vdash i \notin \dot{A}^\alpha$ iff $k^\alpha_{n,i} = 0$, we have $\{p^\alpha_{n,i} \mid i, n \in \omega\} \subseteq \mathbb{P}\restriction_{B^\alpha}$. We may also assume all $B^\alpha$ are trees. Letting $B := \bigcup_{\alpha<\kappa} B^\alpha$, we see that $|B| < \lambda$. By CH and the $\Delta$-system lemma, we may also assume that $\{B^\alpha \mid \alpha < \omega_2\}$ is a $\Delta$-system with root $R$ and that there is a tree $T \subseteq (\omega_1^* \cup \omega_1)^{<\omega}$ such that

- any $B^\alpha$ is weakly isomorphic to $T$,

[4] For an alternative explanation, see [Bre02, Lemma 4.6].

- $\varphi_{\alpha,\beta}\colon B^\alpha \to B^\beta$ is an isomorphism of trees (as defined above) fixing $R$,
- the induced isomorphism $\psi_{\alpha,\beta}\colon \mathbb{P}\restriction_{B^\alpha} \to \mathbb{P}\restriction_{B^\beta}$ maps $p^\alpha_{n,i}$ to $p^\beta_{n,i}$,
- there are numbers $k_{n,i}$ such that $k^\alpha_{n,i} = k_{n,i}$ for all $\alpha < \omega_2$, and
- there is some $\theta_0 < \omega_1$ such that whenever $\alpha < \omega_2$, $x \in B^\alpha$, $j$ odd and $x(j) \in \lambda^*$, then $x(j) \in S^\theta$ for some $\theta < \theta_0$.

By rigidity, we then have $\varphi_{\beta,\alpha} = \varphi^{-1}_{\alpha,\beta}$ and $\varphi_{\alpha,\gamma} = \varphi_{\beta,\gamma} \circ \varphi_{\alpha,\beta}$, and similarly for the $\psi_{\alpha,\beta}$. Further note that the third and fourth clauses immediately imply that $\psi_{\alpha,\beta}$ also maps the name $\dot{A}^\alpha$ to $\dot{A}^\beta$.

For $\alpha < \omega_2$, write $B^\alpha = \{x^\alpha_s \mid s \in T\}$. This means in particular that $|s| = |x^\alpha_s|$, that $s(n)$ is positive iff $x^\alpha_s(n)$ is positive, and that $\varphi_{\alpha,\beta}(x^\alpha_s) = x^\beta_s$. Let $S \subseteq T$ be the subtree corresponding to the root $R$, that is, $s \in S$ iff $x^\alpha_s \in R$ for any $\alpha < \omega_2$. So, for $\alpha \neq \beta$, $x^\alpha_s = x^\beta_s$ iff $s \in S$. List the immediate successors of $S$ in $T$ as $\{t_n \mid n \geq 1\}$, i. e. $\{t_n \mid n \geq 1\} = \{t \in T \smallsetminus S \mid t\restriction_{|t|-1} \in S\}$. For $\alpha < \beta < \omega_2$, define

$$F(\{\alpha,\beta\}) = \begin{cases} n & \text{if } \operatorname{abs}(x^\alpha_{t_n}(|t_n|-1)) > \operatorname{abs}(x^\beta_{t_n}(|t_n|-1)) \\ & \text{(if such an } n \text{ exists and is minimal with this property),} \\ 0 & \text{otherwise.} \end{cases}$$

Note that, by well-foundedness of the ordinals, for every $n \geq 1$, any subset of $\omega_2$ homogeneous in colour $n$ must be finite. Hence, by the Erdős-Rado theorem, we obtain a subset of size $\omega_1$ homogeneous in colour 0 and may as well assume that $\omega_1$ itself is 0-homogeneous.[5] Using further pruning arguments, we may additionally suppose that if $s \in S$ and $(\zeta,\xi) \in (\omega_1^*)^2 \cup (\omega_1)^2$ with $\operatorname{abs}(\zeta) < \operatorname{abs}(\xi)$ and $s^\frown\zeta, s^\frown\xi \in T \smallsetminus S$ (so $s^\frown\zeta = t_n$ and $s^\frown\xi = t_m$ for some $n \neq m \geq 1$), then for all $\alpha < \beta < \omega_1$,

$(\star)$ $\operatorname{abs}(x^\alpha_{s^\frown\zeta}(|s|)) < \operatorname{abs}(x^\beta_{s^\frown\zeta}(|s|))$, all $\operatorname{abs}(x^\alpha_{s^\frown\zeta}(|s|))$ are larger than $\omega_1$, and
- either $\operatorname{abs}(x^\beta_{s^\frown\zeta}(|s|)) < \operatorname{abs}(x^\alpha_{s^\frown\xi}(|s|))$ (which is the case when $\sup_{\alpha<\omega_1} \operatorname{abs}(x^\alpha_{s^\frown\zeta}(|s|)) < \sup_{\alpha<\omega_1} \operatorname{abs}(x^\alpha_{s^\frown\xi}(|s|))$),
- or $\operatorname{abs}(x^\alpha_{s^\frown\xi}(|s|)) < \operatorname{abs}(x^\beta_{s^\frown\zeta}(|s|))$ (which is the case when $\sup_{\alpha<\omega_1} \operatorname{abs}(x^\alpha_{s^\frown\zeta}(|s|)) = \sup_{\alpha<\omega_1} \operatorname{abs}(x^\alpha_{s^\frown\xi}(|s|))$).

We can additionally assume that for $s \in S$ and $(\zeta,\xi) \in (\omega_1^*)^2 \cup (\omega_1)^2$ with $\operatorname{abs}(\zeta) < \operatorname{abs}(\xi)$ and $s^\frown\zeta \in T, s^\frown\xi \in T \smallsetminus S$, if either $s^\frown\zeta \in S$ or the first case of $(\star)$ holds, then

$(\dagger)$ $\sup_{\alpha<\omega_1} \operatorname{abs}(x^\alpha_{s^\frown\zeta}(|s|)) < \varepsilon < \operatorname{abs}(x^0_{s^\frown\xi}(|s|))$ for some $\varepsilon = \varepsilon_{s^\frown\xi}$.

(Note that this $\varepsilon$ depends on $s$ and $\xi$, but not on $\zeta$.) Define $x^\kappa_s \in L$ by recursion on the length of $s \in T$, as follows: If $s \in S$, then let $x^\kappa_s := x^\alpha_s$ for any $\alpha < \omega_1$ (in particular, $|x^\kappa_s| = |x^\alpha_s| = |s|$). If $s \in S$ and $s^\frown\zeta \notin S$, we will have $|x^\kappa_{s^\frown\zeta}| = |s^\frown\zeta| + 2$. First, let $x^\kappa_{s^\frown\zeta}(|s|)$ be the limit of the $x^\alpha_{s^\frown\zeta}(|s|)$ (so it is either the sup or the inf, depending on whether $\zeta$ is positive or negative). Next, find $\gamma < \lambda$ with $\gamma > \omega_1$ and $\gamma^* \in S^{\theta_0}$ such that for all $s$ and $\zeta$, and all $y \in B$ with $y\restriction_{|s|+1} = x^\kappa_{s^\frown\zeta}\restriction_{|s|+1}$, we have $\operatorname{abs}(y(|s|+1)) < \gamma$. It is clear that such a $\gamma$ exists because $\lambda > |B|$ is regular.

$(\star\star)$ If $\zeta$ (and $x^\kappa_{s^\frown\zeta}(|s|)$) is positive, let $x^\kappa_{s^\frown\zeta}(|s|+1) := \gamma^*$, and if $\zeta$ is negative, let $x^\kappa_{s^\frown\zeta}(|s|+1) := \gamma$.

---

[5] Erdős-Rado in the form $(2^\omega)^+ \to ((2^\omega)^+, (\omega_1)_\omega)^2$ actually gives a 0-homogeneous set of size $\omega_2$ and this is needed to guarantee that we can assume the $\operatorname{abs}(x^\alpha_{s^\frown\zeta}(|s|))$ to be larger than $\omega_1$ in $(\star)$ below, but in the end we will use only the first $\omega_1$ many.

To complete the definition of $x^\kappa_{s^\frown\zeta}$, define

$$x^\kappa_{s^\frown\zeta}(|s|+2) = \begin{cases} x^0_{s^\frown\zeta}(|s|) & \text{if } |s|>0, \\ \xi+2n+1 & \text{if } |s|=0 \text{ and } \zeta=\xi+n \text{ with } \xi \text{ limit or } \xi=0. \end{cases}$$

Finally, for the remaining $t\in T$, stipulate again that $|x^\kappa_t|=|t|+2$, find $s\subsetneq t$ with $s\in S$ maximal, let $x^\kappa_t\restriction_{|s|+3} := x^\kappa_{s^\frown t(|s|)}$ and $x^\kappa_t(j+2) := x^0_t(j)$ for $j>|s|$.

Let $B^\kappa := \{x^\kappa_s \mid s\in T\}$. Notice that $B^\kappa$, though very tree-like, is not a tree like the $B^\alpha$. For $\alpha<\omega_1$, define $\varphi_{\alpha,\kappa}\colon B^\alpha\to B^\kappa$ by $\varphi_{\alpha,\kappa}(x^\alpha_s) := x^\kappa_s$ for $s\in T$. We proceed to show that $\mathbb{P}\restriction_{B^\alpha}$ and $\mathbb{P}\restriction_{B^\kappa}$ are isomorphic by the map $\psi_{\alpha,\kappa}$ induced by $\varphi_{\alpha,\kappa}$, which almost maps $\mathcal{I}\restriction_{B^\alpha}$ to $\mathcal{I}\restriction_{B^\kappa}$ (in the sense explained below). It suffices to consider the case $\alpha=0$ as $\mathbb{P}\restriction_{B^\alpha}\cong\mathbb{P}\restriction_{B^0}$. Clearly, $\varphi=\varphi_{0,\kappa}$ is order-preserving, and it is sufficient to figure out the effect of $\varphi$ and its inverse on the $L_\alpha$ and the $J_x$.

First fix $\beta$ and consider $L_\beta$. We claim that there is $\beta_0\le\beta$ such that $\varphi(L_{\beta_0}\cap B^0) = L_{\beta_0}\cap B^\kappa = L_\beta\cap B^\kappa$. To see this, let $\beta'\ge\beta$ be minimal with $\beta'\in\{x^\kappa_s(0) : s\in T\}$ (if there is no such $\beta'$ the claim follows trivially). If $\beta'=x^\kappa_s(0)$ for some $s\in S$ with $|s|=1$ and there is no $t\in T\setminus S$ with $|t|=1$ and $x^\kappa_t(0)=\beta'$, then $x^0_s(0)=\beta'$ as well, and it is easy to see that $\varphi(L_\beta\cap B^0)=L_\beta\cap B^\kappa$ so that $\beta_0=\beta$ works.[6] If $\beta'=x^\kappa_s(0)$ for some $s\in T\setminus S$ with $|s|=1$, then $x^0_s(0)<\beta'$ for any such $s$. In case $\beta'=\beta$, we see by ($\star\star$) that $\varphi(L_\beta\cap B^0)=L_\beta\cap B^\kappa$ again holds so that we can take $\beta_0=\beta$. So assume $\beta<\beta'$. Then the existence of $\beta_0$ (below all such $x^0_s(0)$) is guaranteed by (†) (with $\varepsilon=\varepsilon_{\langle\xi\rangle}$ where $\xi$ is minimal such that $\langle\xi\rangle\in T\setminus S$ and $x^\kappa_{\langle\xi\rangle}(0)=\beta'$).

Now consider the case when $\varphi(L_{\beta_0}\cap B^0)=L_{\beta_0}\cap B^\kappa=L_\beta\cap B^\kappa$, yet $\varphi(L_\beta\cap B^0)\not\subseteq L_\beta$. For any $s\in T$ with $x^0_s\in L_\beta$, but $x^\kappa_s\notin L_\beta$, we must have $x^\kappa_s(0)>\beta\ge x^0_s(0)\ge\beta_0$ and $x^\kappa_s(0)=\sup_{\alpha<\omega_1}x^\alpha_s(0)$. In particular, for all such $s$, $x^\kappa_s(0)$ must have the same value, namely $\beta'$ from the previous paragraph. Moreover, $x^\kappa_s(1)=\gamma^*$ and $x^\kappa_s(2)=\xi+2n+1<\omega_1$, where $s(0)=\xi+n$ with $\xi$ limit or 0. If, for some $s\in T$, $x^0_s(0)=\beta$, let $\eta=\xi+2n+1$, where $s(0)=\xi+n$ with $\xi$ limit or 0. If there is no such $s$ and $\xi+n=\sup\{s(0)+1\mid x^0_s(0)<\beta\}$, $\xi$ limit or 0, let $\eta=\xi+2n$. Then we see that $L_\beta\cap B^0$ is mapped to $(L_\beta\cup J_x)\cap B^\kappa$ via $\varphi$, where $|x|=3$, $x(0)=\beta'$, $x(1)=\gamma^*$, and $x(2)=\eta$ (note that this $x$ is indeed relevant).

Next assume $x$ is relevant and consider $J_x$. Assume that $J_x\cap B^0\neq\varnothing$. Then there must be $s\in T$ such that $|s|=|x|-1$ and $x^0_s=x\restriction_{|x|-1}$. In case $s\in S$, we have $x^\kappa_s=x^0_s$ and $J_x\cap B^0$ is mapped to $J_x\cap B^\kappa$ via $\varphi$ because, by ($\star$), we must have $y\in R$ for any $y\in B^0$ with $|y|=|x|$, $y\restriction_{|x|-1}=x^0_s$ and $y(|x|-1)\le x(|x|-1)<\omega_1$. In case $s\in T\setminus S$, let $j_0<|s|$ be maximal with $s\restriction_{j_0}\in S$. Define $y$ by $|y|:=|x|+2$, $y\restriction_{|y|-1}=x^\kappa_s$ and $y(|y|-1)=x(|x|-1)$ and note that $J_x\cap B^0$ gets mapped to $J_y\cap B^\kappa$ via $\varphi$ provided we can show that $y$ is relevant. In case $j_0>0$, this follows because, whenever $x^0_s(j)>\omega_1$ where $j\ge j_0$ is even, then also $x^\kappa_s(j+2)=x^0_s(j)>\omega_1$, and, if $j_0$ is even, we additionally have $x^\kappa_s(j_0)=\sup_{\alpha<\omega_1}x^\alpha_s(j_0)>\omega_1$ while, if $j_0$ is odd, we additionally have $x^\kappa_s(j_0+1)=\gamma>\omega_1$. In case $j_0=0$, this is true because $x^\kappa_s(1)\in S^{\theta_0}$ and $\theta_0$ is larger than all the $\theta$ for which $x^\kappa_s(j)\in S^\theta$ where $j>1$ is odd.

[6] As pointed out to us by the referee, it is possible that there are $s\in S$ and $t\in T\setminus S$ with $|s|=|t|=1$ and $x^\kappa_s(0)=x^\kappa_t(0)=\beta'$.

On the other hand, assume $J_x \cap B^\kappa \neq \varnothing$. If there is $s \in S$ with $|s| = |x| - 1$ and $x^\kappa_s = x{\restriction}_{|x|-1}$, we conclude as in the previous paragraph that $\varphi$ maps $J_x \cap B^0$ to $J_x \cap B^\kappa$. So suppose there is no such $s \in S$. First, assume $|x| \geq 5$, let $s \in T \setminus S$ with $|s| = |x| - 3$ and $x^\kappa_s = x{\restriction}_{|x|-1}$, and let $j_0 < |s|$ maximal with $s{\restriction}_{j_0} \in S$. Let $y$ be such that $|y| = |x| - 2$, $y{\restriction}_{|y|-1} = x^0_s$, and $y(|y| - 1) = x(|x| - 1)$, and check that $\varphi$ maps $J_y \cap B^0$ to $J_x \cap B^\kappa$ as in the previous paragraph. Finally, assume $|x| = 3$. Then there is a $\beta < x(0)$ such that for any $s \in T \setminus S$ with $x^\kappa_s \in J_x \cap B^\kappa$, we have $x^0_s(0) < \beta$, and we see that $\varphi^{-1}$ maps $J_x \cap B^\kappa$ into $L_\beta \cap B^0$. This is the only case where the templates $\mathcal{I}{\restriction}_{B^0}$ and $\mathcal{I}{\restriction}_{B^\kappa}$ are not identified via $\varphi$, for $\varphi^{-1}(J_x \cap B^\kappa)$ need not belong to $\mathcal{I}{\restriction}_{B^0}$. However, note that only big sets in the template matter for the definition of the iteration, and since $\varphi^{-1}(J_x \cap B^\kappa) \subseteq L_\beta \cap B^0$, it is easy to see that we can conclude that $\mathbb{P}{\restriction}_{B^0} \cong \mathbb{P}{\restriction}_{B^\kappa}$, as witnessed by $\psi_{0,\kappa}$. (For a more formal argument, see [Bre02, Lemma 1.7].)

As mentioned already, this means that $\psi_{\alpha,\kappa}$ is an isomorphism between $\mathbb{P}{\restriction}_{B^\alpha}$ and $\mathbb{P}{\restriction}_{B^\kappa}$, and we can define $\dot{A}^\kappa$ as the $\psi_{\alpha,\kappa}$-image of $\dot{A}^\alpha$ (where $\alpha < \omega_1$ is arbitrary). More explicitly, $p^\kappa_{n,i} = \psi_{\alpha,\kappa}(p^\alpha_{n,i})$, and $p^\kappa_{n,i} \Vdash i \in \dot{A}^\kappa$ iff $k_{n,i} = 1$ and $p^\kappa_{n,i} \Vdash i \notin \dot{A}^\kappa$ iff $k_{n,i} = 0$.

Now assume $F \subseteq \kappa$ is finite. By construction (this is straightforward by our pruning arguments, noting in particular that for all but countably many $\alpha < \omega_1$, we have $B^\alpha \cap \bigcup_{\beta \in F} B^\beta \subseteq R$ and using ($\star\star$)), we see that there is an $\alpha < \omega_1$ such that $B^\alpha \cup \bigcup_{\beta \in F} B^\beta$ and $B^\kappa \cup \bigcup_{\beta \in F} B^\beta$ are order-isomorphic via the mapping $\varphi'$ fixing nodes of $\bigcup_{\beta \in F} B^\beta$ and sending the $x^\alpha_s$ to the corresponding $x^\kappa_s$ via $\varphi_{\alpha,\kappa}$ (in fact, this is also true for all but countably many $\alpha$). Also, by the properties of $\varphi_{\alpha,\kappa}$ explained in the previous paragraphs and the fact that $\varphi'$ is the identity outside of the domain of $\varphi_{\alpha,\kappa}$, $\varphi'$ almost maps $\mathcal{I}{\restriction}_{B^\alpha \cup \bigcup_{\beta \in F} B^\beta}$ to $\mathcal{I}{\restriction}_{B^\kappa \cup \bigcup_{\beta \in F} B^\beta}$ in the sense explained above. Thus the induced map $\psi' \colon \mathbb{P}{\restriction}_{B^\alpha \cup \bigcup_{\beta \in F} B^\beta} \to \mathbb{P}{\restriction}_{B^\kappa \cup \bigcup_{\beta \in F} B^\beta}$ is an isomorphism fixing the names $\dot{A}^\beta$ for $\beta \in F$ and mapping the name $\dot{A}^\alpha$ to the name $\dot{A}^\kappa$. Since $\mathbb{P}$ forces that $\{\dot{A}^\beta \mid \beta \in \{\alpha\} \cup F\}$ is $1/2$-independent, this is actually forced by $\mathbb{P}{\restriction}_{B^\alpha \cup \bigcup_{\beta \in F} B^\beta}$ by complete embeddability. By isomorphism, $\mathbb{P}{\restriction}_{B^\kappa \cup \bigcup_{\beta \in F} B^\beta}$ forces that $\{\dot{A}^\beta \mid \beta \in \{\kappa\} \cup F\}$ is $1/2$-independent, and therefore so does $\mathbb{P}$. Since this holds for every finite $F \subseteq \kappa$, $\mathbb{P}$ actually forces that $\{\dot{A}^\beta \mid \beta \in \kappa + 1\}$ is $1/2$-independent, and thus that $\dot{\mathcal{A}}$ is not maximal, and thus the proof is complete.

We note that this last paragraph is the main difference from the original proof in [Bre21, Theorem 30]. □

We remark that the construction in [Bre03] can be modified analogously to show that $\mathfrak{i}_{1/2}$ can have countable cofinality; see Theorem 4.15 in the subsequent section.

**Corollary 3.13.** Con($\mathfrak{r}_* < \mathfrak{i}_*$).

*Proof.* Replacing the name for a $1/2$-independent family $\dot{\mathcal{A}}$ with the name for a $*$-independent family, the same proof as in Lemma 3.9 shows the analogous result. □

For the final proof of this section, we will require two combinatorial lemmata.

**Lemma 3.14.** *If* $R, S \subseteq \omega$ *are disjoint finite sets of sizes* $r$ *and* $s$*, respectively,* $s = c \cdot r$ *for some* $c > 1$*,* $q \in (0,1)$ *and* $A \subseteq R$*,* $B \subseteq S$ *such that*

$$\frac{|B|}{|S|} \in (q - \varepsilon, q + \varepsilon)$$

*for some* $\varepsilon > 0$*, then*

$$\frac{|A \cup B|}{|R \cup S|} \in \left( q - \varepsilon - \frac{1}{c}, q + \varepsilon + \frac{1}{c} \right).$$

*Proof.* Since

$$\frac{1}{1 + {}^1\!/_c} \geq 1 - \frac{1}{c},$$

we have the lower bound

$$\begin{aligned} \frac{|A \cup B|}{|R \cup S|} &> \frac{s \cdot (q - \varepsilon)}{r + s} = \frac{s \cdot (q - \varepsilon)}{s \cdot {}^1\!/_c + s} = \frac{q - \varepsilon}{1 + {}^1\!/_c} \\ &\geq (q - \varepsilon) \left( 1 - \frac{1}{c} \right) \geq q - \varepsilon - \frac{1}{c}. \end{aligned}$$

For the upper bound, we get

$$\begin{aligned} \frac{|A \cup B|}{|R \cup S|} &< \frac{r + s \cdot (q + \varepsilon)}{r + s} = \frac{s \cdot {}^1\!/_c + s \cdot (q + \varepsilon)}{s \cdot {}^1\!/_c + s} \\ &= \frac{q + \varepsilon + {}^1\!/_c}{1 + {}^1\!/_c} \leq q + \varepsilon + \frac{1}{c}. \end{aligned}$$

□

**Lemma 3.15.** *If* $R, S \subseteq \omega$*,* $0 < r < 1$*,* $\varepsilon > 0$ *and* $m \leq \ell$ *are such that*

$$\frac{|R \cap m|}{m}, \frac{|S \cap m|}{m}, \frac{|S \cap \ell|}{\ell} \in (r - \varepsilon, r + \varepsilon)$$

*then*

$$\frac{|(R \cap m) \cup (S \cap [m, \ell))|}{\ell} \in (r - 3\varepsilon, r + 3\varepsilon).$$

*Proof.* Suppose this were false for some $\ell^* \geq m$; then without loss of generality,

$$\frac{|(R \cap m) \cup (S \cap [m, \ell^*))|}{\ell^*} \geq r + 3\varepsilon.$$

Since

$$\frac{|R \cap m|}{m} < r + \varepsilon,$$

we get

$$\frac{|S \cap [m, \ell^*)|}{\ell^*} \geq r + 3\varepsilon - \frac{m}{\ell^*}(r + \varepsilon).$$

But then

$$\frac{|S \cap m|}{m} > r - \varepsilon$$

implies

$$\frac{|S \cap \ell^*|}{\ell^*} = \frac{|(S \cap m) \cup (S \cap [m, \ell^*))|}{\ell^*} > \frac{m}{\ell^*}(r - \varepsilon) + r + 3\varepsilon - \frac{m}{\ell^*}(r + \varepsilon)$$
$$= r + 3\varepsilon - \frac{2m}{\ell^*} \cdot \varepsilon \geq r + \varepsilon,$$

which is a contradiction. □

**Lemma 3.16.** $\mathrm{Con}(\mathfrak{i}_{1/2} < \mathrm{cov}(\mathcal{M}))$ *and thus* $\mathrm{Con}(\mathfrak{i}_{1/2} < \mathfrak{i})$.

*Proof.* The proof is analogous to the classical proof of $\mathrm{Con}(\aleph_1 = \mathfrak{a} < 2^{\aleph_0})$ (see e. g. [Hal17, Proposition 18.5]).

Assume CH in the ground model and let $\lambda \geq \aleph_2$ with $\lambda = \lambda^{\aleph_0}$. We force with the $\lambda$-Cohen forcing poset $\mathbb{C}_\lambda$; letting $G$ be a $\mathbb{C}_\lambda$-generic filter, it is clear that $V[G] \vDash \mathrm{cov}(\mathcal{M}) = \mathfrak{i} = 2^{\aleph_0} = \lambda$. We will now show $V[G] \vDash \mathfrak{i}_{1/2} = \aleph_1$ by constructing a maximal $1/2$-independent family $\mathcal{A}$ in the ground model such that $\mathcal{A}$ remains maximal $1/2$-independent in $V[G]$. By the usual arguments, it suffices to consider what happens to a countably infinite $1/2$-independent family when forcing with just $\mathbb{C} := \langle 2^{<\omega}, \subseteq \rangle$.

Let $\mathcal{A}_0 := \{A_n \subseteq [\omega]^{\aleph_0} \mid n < \omega\}$ be such a family. Fix (in the ground model) an enumeration $\{(p_\alpha, \dot{X}_\alpha) \mid \omega \leq \alpha < \omega_1\}$ of all pairs $(p, \dot{X})$ such that $p \in \mathbb{C}$ and $\dot{X}$ is a nice name for a subset of $\omega$.[7] In particular, this means that for any $\langle \check{n}, p_1 \rangle, \langle \check{n}, p_2 \rangle \in \dot{X}$, either $p_1 = p_2$ or $p_1 \perp p_2$. Note that since $V \vDash$ CH, there are just $\aleph_1$ many nice names for subsets of $\omega$ in $V$.

We now construct $\mathcal{A}$ from $\mathcal{A}_0$ iteratively as follows: Let $\omega \leq \alpha < \omega_1$ and assume we have already defined sets $A_\beta \subseteq \omega$ for all $\beta < \alpha$. Below, we will construct $A_\alpha \subseteq \omega$ such that the following two properties hold:

(i) The family $\{A_\beta \mid \beta \leq \alpha\}$ is $1/2$-independent.
(ii) If $p_\alpha \Vdash |\dot{X}_\alpha| = \aleph_0 \wedge$ "$\{A_\beta \mid \beta < \alpha\} \cup \{\dot{X}_\alpha\}$ is $1/2$-independent", then for all $m < \omega$, the set $D_m^\alpha := \{q \in \mathbb{C} \mid \exists n \geq m \colon q \Vdash A_\alpha \cap [2^n, 2^{n+1}) = \dot{X}_\alpha \cap [2^n, 2^{n+1})\}$ is dense below $p_\alpha$.

We first show that the $\mathcal{A} := \{A_\beta \mid \beta < \omega_1\}$ constructed this way is a maximal $1/2$-independent family in $V^{\mathbb{C}}$. Clearly, $\mathcal{A}$ is $1/2$-independent, so only maximality could fail. Suppose it were not maximal; then there is a condition $p$ and a nice name $\dot{X}$ for a subset of $\omega$ such that $p \Vdash$ "$\mathcal{A} \cup \{\dot{X}\}$ is $1/2$-independent". Let $\alpha$ be such that $(p, \dot{X}) = (p_\alpha, \dot{X}_\alpha)$ and let $\varepsilon > 0$ be sufficiently small (e. g. $\varepsilon < 1/16$). We can then find $q \leq p_\alpha$ and $m < \omega$ such that

$$(*_1) \qquad q \Vdash \frac{|A_\alpha \cap \dot{X}_\alpha \cap \ell|}{\ell} \in \left(\frac{1}{4} - \varepsilon, \frac{1}{4} + \varepsilon\right) \text{ for all } \ell \geq 2^m$$

(because $p_\alpha$ forces that $\{A_\alpha, \dot{X}_\alpha\}$ is $1/2$-independent) and

$$\frac{|A_\alpha \cap [2^n, 2^{n+1})|}{2^n} > \frac{1}{2} - \varepsilon \text{ for all } n \geq m.$$

[7] The reason the index set of the enumeration is $[\omega, \omega_1)$ instead of $[0, \omega_1)$ is just to make the notation more convenient.

Now by the density of $D^\alpha_m$ below $p_\alpha$, we can find $r \leq q$ and some $n \geq m$ such that $r \Vdash A_\alpha \cap [2^n, 2^{n+1}) = \dot{X}_\alpha \cap [2^n, 2^{n+1})$. But this implies that

$$r \Vdash \frac{|A_\alpha \cap \dot{X}_\alpha \cap 2^{n+1}|}{2^{n+1}} = \frac{1}{2} \cdot \frac{|A_\alpha \cap \dot{X}_\alpha \cap 2^n|}{2^n} + \frac{1}{2} \cdot \frac{|A_\alpha \cap \dot{X}_\alpha \cap [2^n, 2^{n+1})|}{2^n}$$
$$> \frac{1/4 - \varepsilon}{2} + \frac{1/2 - \varepsilon}{2} = \frac{3}{8} - \varepsilon > \frac{1}{4} + \varepsilon,$$

which contradicts Eq. ($*_1$).

We finally have to show that we can find such an $A_\alpha$ satisfying (i) and (ii) for any $\omega \leq \alpha < \omega_1$. We only have to consider those $\alpha$ such that $\dot{X}_\alpha$ satisfies the assumption in property (ii), since finding an $A_\alpha$ with property (i) is straightforward. Enumerate $\{A_\beta \mid \beta < \alpha\}$ as $\{B_n \mid n < \omega\}$. For $n < \omega$ and any partial function $f\colon n \to \{-1, 1\}$, we let

$$B^f := \bigcap_{i \in \operatorname{dom}(f)} B_i^{f(i)},$$

where $B^1_i := B_i$ and $B^{-1}_i := \omega \setminus B_i$. We further pick some strictly decreasing sequence of real numbers $\langle \delta_n \mid n < \omega \rangle$ with $\delta_0 := 3$ and $\lim_{n\to\infty} \delta_n = 0$ and let $\langle q_n \mid n < \omega \rangle$ be some sequence enumerating all conditions below $p_\alpha$ infinitely often. We will now construct, by induction on $n < \omega$, conditions $r_n \leq q'_n \leq q_n$, a strictly increasing sequence of natural numbers $\langle k_n \mid n < \omega \rangle$ and initial segments $Z_n = A_\alpha \cap 2^{k_n}$ of $A_\alpha$ such that for all $n < \omega$ and all partial functions $f\colon n \to \{-1, 1\}$, the following four statements will hold (with $F := |\operatorname{dom}(f)| + 1$)

(R1) $\dfrac{|B^f \cap Z_n \cap 2^{k_n}|}{2^{k_n}}, \dfrac{|(B^f \setminus Z_n) \cap 2^{k_n}|}{2^{k_n}} \in \left(\dfrac{1}{2^F} - \dfrac{\delta_n}{3}, \dfrac{1}{2^F} + \dfrac{\delta_n}{3}\right),$

(R2) $q'_n \Vdash \dfrac{|B^f \cap \dot{X}_\alpha \cap \ell|}{\ell}, \dfrac{|(B^f \setminus \dot{X}_\alpha) \cap \ell|}{\ell} \in \left(\dfrac{1}{2^F} - \dfrac{\delta_n}{3}, \dfrac{1}{2^F} + \dfrac{\delta_n}{3}\right)$
for all $\ell$ with $2^{k_n} \leq \ell$,

(R3) $\dfrac{|B^f \cap Z_{n+1} \cap \ell|}{\ell}, \dfrac{|(B^f \setminus Z_{n+1}) \cap \ell|}{\ell} \in \left(\dfrac{1}{2^F} - \delta_n, \dfrac{1}{2^F} + \delta_n\right)$
for all $\ell$ with $2^{k_n} \leq \ell \leq 2^{k_{n+1}}$, and

(R4) $r_n \Vdash Z_{n+1} \cap [2^{k_n}, 2^{k_{n+1}}) = \dot{X}_\alpha \cap [2^{k_n}, 2^{k_{n+1}})$.

It is clear that (R1)–(R4) taken together for all $n < \omega$ imply that $A_\alpha := \bigcup_{n\in\omega} Z_n$ is as required by (i) and (ii).

For $n = 0$, let $k_0 := 0$, $q'_0 := q_0$ and $Z_0 := \varnothing$; then (R1) and (R2) hold vacuously by our choice of $\delta_0$, and there is nothing to show yet for (R3) and (R4).

Now assume that we have obtained $k_n$, $q'_n \leq q_n$ and $Z_n$ such that (R1) and (R2) hold for $n$; we will construct $r_n \leq q'_n$, $k_{n+1}$, $q'_{n+1} \leq q_{n+1}$ and $Z_{n+1}$ such that (R3) and (R4) hold for $n$ and such that (R1) and (R2) hold for $n + 1$. We first find $q'_{n+1} \leq q_{n+1}$ and $k'_n \geq k_n$ such that for all partial functions $f\colon n + 1 \to \{-1, 1\}$, we have that (with $F := |\operatorname{dom}(f)| + 1$)

$$q'_{n+1} \Vdash \frac{|B^f \cap \dot{X}_\alpha \cap \ell|}{\ell}, \frac{|(B^f \setminus \dot{X}_\alpha) \cap \ell|}{\ell} \in \left(\frac{1}{2^F} - \frac{\delta_{n+1}}{3}, \frac{1}{2^F} + \frac{\delta_{n+1}}{3}\right)$$

for all $\ell \geq 2^{k'_n}$ (hence satisfying (R2) for $n+1$); this is possible since the assumption in property (ii) is true. Next we find $r_n \leq q'_n$ and a sufficiently large $k_{n+1} \geq k'_n$ such that for all partial functions $f \colon n+1 \to \{-1, 1\}$, we have that (still with $F := |\operatorname{dom}(f)| + 1$)

$$(*_2) \qquad r_n \Vdash \frac{|B^f \cap \dot{X}_\alpha \cap 2^{k_{n+1}}|}{2^{k_{n+1}}}, \frac{|(B^f \setminus \dot{X}_\alpha) \cap 2^{k_{n+1}}|}{2^{k_{n+1}}} \in \left(\frac{1}{2^F} - \frac{\delta_{n+1}}{6}, \frac{1}{2^F} + \frac{\delta_{n+1}}{6}\right)$$

and that $r_n$ decides $\dot{X}_\alpha \cap 2^{k_{n+1}}$; in particular, let $X_n \subseteq [2^{k_n}, 2^{k_{n+1}})$ be such that $r_n \Vdash \dot{X}_\alpha \cap [2^{k_n}, 2^{k_{n+1}}) = X_n$. All this is also possible since the assumption in property (ii) is true. Let $Z_{n+1} := Z_n \cup X_n$.

Now, (R4) holds for $n$ by definition of $Z_{n+1}$. Let $W$ be such that $r_n \Vdash \dot{X}_\alpha \cap 2^{k_{n+1}} = W$. Apply Lemma 3.15 to $R := B^f \cap Z_n$ or $B^f \setminus Z_n$, $S := B^f \cap W$ or $B^f \setminus W$, $r := 1/2^F$, $\varepsilon := \delta_n/3$, $m := 2^{k_n}$ and any $\ell$ with $2^{k_n} \leq \ell \leq 2^{k_{n+1}}$ to see that (R3) for $n$ follows from (R1) and (R2) for $n$ and our choice of $Z_{n+1}$. Finally, apply Lemma 3.14 to $R := 2^{k_n} = r$, $S := [2^{k_n}, 2^{k_{n+1}})$, $s = 2^{k_{n+1}} - 2^{k_n}$, $c = 2^{k_{n+1}-k_n} - 1$, $q = 1/2^F$ and $\varepsilon =: \delta_{n+1}/6$ to see that (R1) for $n+1$ follows from Eq. ($*_2$), (R4) for $n$ and the choice of $k_{n+1}$ sufficiently large as to guarantee $1/c < \delta_{n+1}/6$.

By the usual arguments, our construction implies that $\mathcal{A}$ remains maximal $1/2$-independent in $V^{\mathbb{C}_\lambda}$. $\square$

## 4. More on $\mathfrak{i}_{1/2}$

We describe a forcing for adding a maximal $1/2$-independent family generically with a product-style forcing (like Hechler's forcing for adding a mad family [Hec72]). This gives an alternative proof of the consistency of $\mathfrak{i}_{1/2} < \mathfrak{c}$, while also showing that there can be (consistently) simultaneously maximal $1/2$-independent families of many different sizes and that $\operatorname{cf}(\mathfrak{i}_{1/2}) = \omega$ is consistent. We note in this context that the consistency of $\operatorname{cf}(\mathfrak{i}) = \omega$ is a well-known open problem.

**Definition 4.1.** Fix an uncountable cardinal $\kappa$. We define the forcing $\mathbb{P} = \mathbb{P}_\kappa$ as follows. Conditions are of the form $p = (F^p, n^p, \bar{a}^p, \varepsilon^p)$ such that

(C1) $F^p \subseteq \kappa$ is finite,
(C2) $n^p \in \omega$,
(C3) $\bar{a}^p = \langle a^p_\alpha \subseteq n^p \mid \alpha \in F^p \rangle$,
(C4) $\varepsilon^p \colon 2^{\leq F^p} \longrightarrow \mathbb{Q}^+$ (where $2^{\leq F^p}$ denotes the partial functions from $F^p$ to 2) is such that $\varepsilon^p(f) \leq \varepsilon^p(g)$ whenever $f \subseteq g$,
(C5) for all $f \in 2^{\leq F^p}$, we have

$$\left| \frac{|\bigcap_{f(\alpha)=1} a^p_\alpha \cap \bigcap_{f(\alpha)=0} (n^p \setminus a^p_\alpha)|}{n^p} - \frac{1}{2^{|\operatorname{dom}(f)|}} \right| < \frac{\varepsilon^p(f)}{8},$$

and

(C6) we have

$$\frac{2^{2|F^p|}}{n^p} < \frac{\varepsilon^p}{8},$$

where $\varepsilon^p := \varepsilon^p(\varnothing) = \min\{\varepsilon^p(f) \mid f \in 2^{\leq F^p}\}$

The order is given by $q \leq p$ if

(D1) $F^p \subseteq F^q$,
(D2) $n^p \leq n^q$,
(D3) $a^p_\alpha = a^q_\alpha \cap n^p$ for all $\alpha \in F^p$,
(D4) $\varepsilon^p(f) \geq \varepsilon^q(f)$ for all $f \in 2^{\leq F^p}$, and
(D5) for all $i$ with $n^p \leq i \leq n^q$ and all $f \in 2^{\leq F^p}$, we have

$$\left| \frac{|\bigcap_{f(\alpha)=1}(i \cap a^q_\alpha) \cap \bigcap_{f(\alpha)=0}(i \setminus a^q_\alpha)|}{i} - \frac{1}{2^{|\mathrm{dom}(f)|}} \right| < \varepsilon^p(f).$$

We first need to check we can extend conditions arbitrarily.

**Definition 4.2.** Given a condition $p$ and $E \subseteq \kappa$, we define the restriction $p' = p{\restriction}_E$ of $p$ to $E$ by

(i) $F^{p'} = F^p \cap E$,
(ii) $n^{p'} = n^p$,
(iii) $a^{p'}_\alpha = a^p_\alpha$ for $\alpha \in F^{p'}$, and
(iv) $\varepsilon^{p'} = \varepsilon^p{\restriction}_{2^{\leq F^{p'}}}$.

It is easy to see that $p' \in \mathbb{P}$ and that $p \leq p'$. Also, for $f \in 2^{\leq F^p}$, let

$$b^p_f := \bigcap_{f(\alpha)=1} a^p_\alpha \cap \bigcap_{f(\alpha)=0} (n^p \setminus a^p_\alpha).$$

**Lemma 4.3** (extendibility lemma). *Let $p \in \mathbb{P}$, $E \subseteq \kappa$, $p' = p{\restriction}_E$, $m \in \omega$, and $\varepsilon\colon 2^{\leq F^p} \longrightarrow \mathbb{Q}^+$ with $\varepsilon(f) \leq \varepsilon(g)$ whenever $f \subseteq g$ and $\varepsilon(f) \leq \varepsilon^p(f)$ for all $f \in 2^{\leq F^p}$. Assume $q' \leq p'$ is such that $F^{q'} \subseteq E$. Then there is a condition $q \in \mathbb{P}$ with $q \leq p$, $q \leq q'$, $F^q = F^p \cup F^{q'}$, $n^q \geq m$, and*

- *$\varepsilon^q(f) = \min\{\varepsilon(f), \varepsilon^{q'}(f)\}$ for all $f \in 2^{\leq F^{p'}}$,*
- *$\varepsilon^q(f) = \varepsilon(f)$ for all $f \in 2^{\leq F^p} \setminus 2^{\leq F^{p'}}$,*
- *$\varepsilon^q(f) = \varepsilon^{q'}(f)$ for all $f \in 2^{\leq F^{q'}} \setminus 2^{\leq F^{p'}}$, and*
- *$\varepsilon^q(f) \geq 16$ for all other $f \in 2^{\leq F^q}$.*

*Proof.* Let $F := F^q := F^{q'} \cup F^p$. Define $\varepsilon^q\colon 2^{\leq F} \longrightarrow \mathbb{Q}^+$ as stipulated in the statement of the lemma. Finally, let $n := n^q \geq \max\{m, n^{q'}\}$ be so large that

- $n - n^{q'}$ is divisible by $2^{|F|}$,
- $\frac{n^{q'}}{n} < \frac{\varepsilon^q}{8}$, and
- $\frac{2^{2|F|}}{n} < \frac{\varepsilon^q}{8}$.

Note that the last item immediately guarantees (C6). We produce the required extension in two steps. The main point is to prove (D5) for $q \leq p$ and $q \leq q'$ and condition (C5) for $q \in \mathbb{P}$.

In the first step we extend to $n^{q'}$. This step is only necessary if $E \neq \varnothing$ and $n^{q'} > n^p$. Let $\{\alpha_\ell \mid \ell \in |F^p \smallsetminus E|\}$ enumerate $F^p \smallsetminus E$. For each $f \in 2^{F^{p'}}$, let $c_f := b^{q'}_f \smallsetminus b^p_f = b^{q'}_f \smallsetminus n^p$. Note that the $c_f$ are pairwise disjoint, that their union is the interval $[n^p, n^{q'})$ and that in case $F^{p'} = \varnothing$, we have $c_\varnothing = [n^p, n^{q'})$.

Let $\{c_f(j) \mid j \in m_f\}$ be the increasing enumeration of $c_f$. For each $\ell \in |F^p \smallsetminus E|$ and each $f \in 2^{F^{p'}}$, define

$$a^q_{\alpha_\ell} \cap c_f := \left\{ c_f(j) \;\middle|\; j \in m_f \cap \bigcup_k [2^{\ell+1}k, 2^{\ell+1}k + 2^\ell) \right\}. \tag{$*_3$}$$

Thus $a^q_{\alpha_\ell} \cap [n^p, n^{q'})$ is the disjoint union of the sets $a^q_{\alpha_\ell} \cap c_f$. We need to see that (D5) is satisfied for all $i$ with $n^p \leq i \leq n^{q'}$ and all $g \in 2^{\leq F^p}$. Hence we fix such $i$ and $g$. We may assume that $\operatorname{dom}(g) \not\subseteq E$ (otherwise, (D5) holds by $q' \leq p'$). We will only show that

$$\frac{|i \cap b^q_g|}{i} < \frac{1}{2^{|\operatorname{dom}(g)|}} + \varepsilon^p(g);$$

the second inequality is analogous.

Let $f = g{\restriction}_E = g{\restriction}_{F^{p'}} \in 2^{\leq F^{p'}}$, hence $f \subsetneq g$. By (C5) for $p$ and $f$, we know that

$$\left| n^p \cap b^{q'}_f \right| = \left| b^p_f \right| > n^p \cdot \left( \frac{1}{2^{|\operatorname{dom}(f)|}} - \frac{\varepsilon^p(f)}{8} \right),$$

and by (D5) for $q' \leq q$ and $f$,

$$\left| i \cap b^{q'}_f \right| < i \cdot \left( \frac{1}{2^{|\operatorname{dom}(f)|}} + \varepsilon^p(f) \right);$$

thus

$$\left| [n^p, i) \cap b^{q'}_f \right| < \frac{i - n^p}{2^{|\operatorname{dom}(f)|}} + \frac{9i \cdot \varepsilon^p(f)}{8}.$$

For $f' \in 2^{F^{p'}}$ with $f \subseteq f'$ we have, by Eq. ($*_3$),

$$\left| [n^p, i) \cap b^q_{f' \cup g} \right| = \left| i \cap c_{f'} \cap b^q_{g{\restriction}_{F^p \smallsetminus E}} \right| \leq \frac{1}{2^{|\operatorname{dom}(g) \smallsetminus E|}} \cdot |i \cap c_{f'}| + 2^{|F^p \smallsetminus E|}.$$

Indeed, if $i$ is such that $i \cap c_{f'} = \{c_{f'}(j) \mid j \in \bar{m}\}$, where $\bar{m} \leq m_f$ is divisible by $2^{|F^p \smallsetminus E|}$, then the partition in Eq. ($*_3$) yields that the set on the left-hand side has size exactly $\frac{1}{2^{|\operatorname{dom}(g) \smallsetminus E|}} \cdot |i \cap c_{f'}|$. Therefore, for other $i$, the error is at most $2^{|F^p \smallsetminus E|}$ and the inequality follows.

Since $[n^p, i) \cap b^q_g$ is the disjoint union of the $[n^p, i) \cap b^q_{f' \cup g}$ and $[n^p, i) \cap b^{q'}_f$ is the disjoint union of the $i \cap c_{f'}$, we see that

$$\begin{aligned}
\left|[n^p, i) \cap b^q_g\right| &= \sum_{f \subseteq f' \in 2^{F^{p'}}} \left|[n^p, i) \cap b^q_{f' \cup g}\right| \\
&\leq \frac{1}{2^{|\mathrm{dom}(g) \smallsetminus E|}} \cdot \sum_{f \subseteq f' \in 2^{F^{p'}}} |i \cap c_{f'}| + 2^{|F^{p'} \smallsetminus \mathrm{dom}(f)|} \cdot 2^{|F^p \smallsetminus E|} \\
&\leq \frac{1}{2^{|\mathrm{dom}(g) \smallsetminus E|}} \cdot \left|[n^p, i) \cap b^{q'}_f\right| + 2^{|F^p|} \\
&< \frac{i - n^p}{2^{|\mathrm{dom}(g)|}} + \frac{9i \cdot \varepsilon^p(f)}{8 \cdot 2^{|\mathrm{dom}(g) \smallsetminus E|}} + 2^{|F^p|}
\end{aligned}$$

and thus, by (C5) for $p$ and $g$ and (C6) for $p$, and using that $g$ strictly extends $f$,

$$(*_4) \qquad \begin{aligned}
\frac{|i \cap b^q_g|}{i} &= \frac{|n^p \cap b^q_g|}{i} + \frac{|[n^p, i) \cap b^q_g|}{i} \\
&< \frac{1}{2^{|\mathrm{dom}(g)|}} + \frac{\varepsilon^p(g)}{8} + \frac{9 \cdot \varepsilon^p(f)}{16} + \frac{\varepsilon^p}{8} < \frac{1}{2^{|\mathrm{dom}(g)|}} + \frac{7 \cdot \varepsilon^p(g)}{8},
\end{aligned}$$

as required.

We now extend from $n^{q'}$ to $n = n^q$. Let $\{\alpha_\ell \mid \ell \in |F^{p'}|\}$ enumerate $F^{p'}$. Next let $\tilde{\ell} = \min\{|F^p \smallsetminus F^{p'}|, |F^{q'} \smallsetminus F^{p'}|\}$. Let $\{\alpha_{2\ell + |F^{p'}|} \mid \ell < \tilde{\ell}\}$ enumerate the next $\tilde{\ell}$ many elements of $F^p \smallsetminus F^{p'} = F^p \smallsetminus E$, and let $\{\alpha_{2\ell+1+|F^{p'}|} \mid \ell < \tilde{\ell}\}$ enumerate the next $\tilde{\ell}$ many elements of $F^{q'} \smallsetminus F^{p'}$. Finally let $\{\alpha_\ell \mid |F^{p'}| + 2\tilde{\ell} \leq \ell < |F|\}$ enumerate the remaining elements of $F$. Define

$$(*_5) \qquad a^q_{\alpha_\ell} \cap [n^{q'}, n) = \bigcup_k [n^{q'} + 2^{\ell+1}k, n^{q'} + 2^{\ell+1}k + 2^\ell)$$

for $\ell < |F|$. First, we need to show (D5) for all $i$ with $n^{q'} \leq i < n$ and all $g \in 2^{\leq F^p} \cup 2^{\leq F^{q'}}$. Fix such $i$ and $g$. Without loss of generality, we may assume $g \in 2^{\leq F^p}$. (For $g \in 2^{\leq F^{q'}}$ the proof is the same.) Again, we only show the inequality

$$\frac{|i \cap b^q_g|}{i} < \frac{1}{2^{|\mathrm{dom}(g)|}} + \varepsilon^p(g).$$

By Eq. ($*_5$) and the choice of the sequence of the $\alpha_\ell$, we have

$$\left|[n^{q'}, i) \cap b^q_g\right| \leq \frac{i - n^{q'}}{2^{|\mathrm{dom}(g)|}} + 2^{2|F^p|}.$$

Thus, by Eq. ($*_4$) for $n^{q'}$, we have

$$\begin{aligned}
\frac{|i \cap b^q_g|}{i} &= \frac{|n^{q'} \cap b^q_g|}{i} + \frac{|[n^{q'}, i) \cap b^q_g|}{i} \\
&< \frac{1}{2^{|\mathrm{dom}(g)|}} + \frac{7 \cdot \varepsilon^p(g)}{8} + \frac{2^{2|F^p|}}{i} < \frac{1}{2^{|\mathrm{dom}(g)|}} + \varepsilon^p(g),
\end{aligned}$$

as required.

Finally, we need to show condition (C5) for $q$ and $g \in 2^{\leq F}$. Since $n - n^{q'}$ is divisible by $2^{|F|}$, it is easy to see that

$$\left| [n^{q'}, n) \cap b^p_g \right| = \frac{n - n^{q'}}{2^{|\mathrm{dom}(g)|}}.$$

Thus

$$\frac{1}{2^{|\mathrm{dom}(g)|}} \cdot \frac{n - n^{q'}}{n} \leq \frac{|b^p_g|}{n} \leq \frac{1}{2^{|\mathrm{dom}(g)|}} \cdot \frac{n - n^{q'}}{n} + \frac{n^{q'}}{n},$$

and the required inequality follows from $\frac{n^{q'}}{n} < \frac{\varepsilon^q}{8}$. □

**Corollary 4.4.** *Let $p \in \mathbb{P}$ and $m \in \omega$. Then there is a condition $q \in \mathbb{P}$ with $q \leq p$ and $n^q \geq m$. Furthermore, we may require $F^q = F^p$ and $\varepsilon^q = \varepsilon^p$.*

*Proof.* Apply Lemma 4.3 with $E = \varnothing$ (so $p' = q'$ is a trivial condition[8]) and $\varepsilon = \varepsilon^p$. □

**Corollary 4.5.** *Let $p \in \mathbb{P}$ and $\alpha \in \kappa$. Then there is a condition $q \in \mathbb{P}$ with $q \leq p$ and $\alpha \in F^q$.*

*Proof.* We may assume $\alpha \notin F^p$. Apply Lemma 4.3 with $E = \{\alpha\}$ (so $p'$ is a trivial condition) and arbitrary $q'$ with $F^{q'} = E = \{\alpha\}$. □

**Corollary 4.6.** *Let $p \in \mathbb{P}$ and $\varepsilon \colon 2^{\leq F^p} \longrightarrow \mathbb{Q}^+$ with $\varepsilon(f) \leq \varepsilon(g)$ whenever $f \subseteq g$. Then there is a condition $q \in \mathbb{P}$ with $q \leq p$ such that $\varepsilon^q(f) \leq \varepsilon(f)$ for all $f \in 2^{\leq F^p}$.*

*Proof.* Apply Lemma 4.3 with $E = \varnothing$ (so $p' = q'$ is a trivial condition). □

**Lemma 4.7** (compatibility lemma)**.** *Assume $p, q \in \mathbb{P}$ are such that $n^p = n^q$, $a^p_\alpha = a^q_\alpha$ for all $\alpha \in F^p \cap F^q$, and $\varepsilon^p \restriction_{2^{\leq (F^p \cap F^p)}} = \varepsilon^q \restriction_{2^{\leq (F^p \cap F^p)}}$. Then $p$ and $q$ are compatible.*

*Proof.* Apply Lemma 4.3 with $p = p$, $E = F^q$, $m = n^p$, and $\varepsilon = \varepsilon^p$. Note that $q' = q$ satisfies the necessary assumptions. □

**Corollary 4.8** (ccc)**.** *$\mathbb{P}$ satisfies Knaster's condition (and thus is ccc) and therefore preserves cardinals.*

*Proof.* This follows from a $\Delta$-system argument together with Lemma 4.7. □

**Definition 4.9.** For $X \subseteq \kappa$, let $\mathbb{P}_X$ be the collection of conditions $p \in \mathbb{P}_\kappa$ with $F^p \subseteq X$.

**Corollary 4.10** (complete embeddability)**.** *For any $X \subseteq \kappa$, $\mathbb{P}_X$ completely embeds into $\mathbb{P}_\kappa$.*

*Proof.* By Lemma 4.3, $p' = p\restriction_X \in \mathbb{P}_X$ is a reduction of $p \in \mathbb{P}_\kappa$. □

[8] The forcing $\mathbb{P}$ does not contain a single trivial condition because there are many conditions with empty $F^p$, but with different $n^p$ and $\varepsilon^p(\varnothing)$. However, all these trivial conditions are identified with the maximal element in the complete Boolean algebra associated with $\mathbb{P}$.

Note that since $\mathbb{P}_\omega$ is countable, it is forcing-equivalent to Cohen forcing $\mathbb{C}_\omega$, and $\mathbb{P}_{\omega_1}$ is forcing-equivalent to the partial order $\mathbb{C}_{\omega_1}$ adding $\omega_1$ many Cohen reals, by Corollary 4.10 and well-known arguments (see e. g. [BJZ97, Theorem 3.2]).[9]

Let $G$ be $\mathbb{P}$-generic over $V$. For $\alpha < \kappa$, let $A_\alpha := \bigcup\{a^p_\alpha \mid p \in G\}$. By the corollaries of Lemma 4.3 (Corollary 4.4, Corollary 4.5 and Corollary 4.6), we immediately see:

**Corollary 4.11.** *$\{A_\alpha \mid \alpha < \kappa\}$ is a $1/2$-independent family.*

Next, combining the basic idea of Hechler's classical work [Hec72] with the combinatorics of Lemma 4.3, we have:

**Lemma 4.12** (maximality)**.** *$\{A_\alpha \mid \alpha < \kappa\}$ is a maximal $1/2$-independent family. Moreover, for any ccc forcing $\mathbb{Q}$ in $V$, $\{A_\alpha \mid \alpha < \kappa\}$ is still maximal in the $\mathbb{P} \times \mathbb{Q}$-generic extension.*

*Proof.* Let $\dot{B}$ be a $\mathbb{P} \times \mathbb{Q}$-name for an infinite and coinfinite subset of $\omega$. For each $i \in \omega$, let $M_i$ be a maximal antichain of conditions deciding $i \in \dot{B}$. By Corollary 4.8, each $M_i$ is at most countable because the product of a ccc forcing and a forcing satisfying Knaster's condition is ccc. Thus we can find a countable $X \subseteq \kappa$ such that $F^p \subseteq X$ for all $(p, \bar{p}) \in \bigcup_i M_i$. Let $\beta \in \kappa \setminus X$. Clearly, it suffices to show:

**Claim.** *Assume $(p_0, \bar{p}_0) \in \mathbb{P} \times \mathbb{Q}$ forces that $\dot{B}$ is $1/2$-independent from all $\dot{A}_\alpha$ for $\alpha \in X$. Then $(p_0, \bar{p}_0)$ forces that for all $k$, there is an $\ell > k$ such that*

$$\frac{|\ell \cap \dot{B} \cap \dot{A}_\beta|}{\ell} > \frac{3}{8}.$$

(Note that, analogously, we can show that $(p_0, \bar{p}_0)$ forces that for all $k$ there is an $\ell > k$ such that

$$\frac{|\ell \cap \dot{B} \cap \dot{A}_\beta|}{\ell} < \frac{1}{8},$$

and in fact, it is not difficult to see that an elaboration of the argument shows that $(p_0, \bar{p}_0)$ forces $\underline{d}(\dot{B} \cap \dot{A}_\beta) = 0$ and $\bar{d}(\dot{B} \cap \dot{A}_\beta) = 1/2$.)

Fix $(p, \bar{p}) \leq (p_0, \bar{p}_0)$ in $\mathbb{P} \times \mathbb{Q}$ and $k$. We need to find $\ell > k$ and $(r, \bar{r}) \leq (p, \bar{p})$ forcing the required statement. We may assume $n^p \geq k$ and $\beta \in F^p$. We may also assume that for $f_0$ with $\mathrm{dom}(f_0) = \{\beta\}$ and $f_0(\beta) = 1$, $\varepsilon^p(f_0) < 1/2$.

Let $p' = p{\restriction}_X$. For $f \in 2^{\leq F^{p'} \cup \{\beta\}}$ with $\beta \in \mathrm{dom}(f)$, let $\dot{C}_f$ denote the $\mathbb{P} \times \mathbb{Q}$-name

$$\bigcap_{\substack{\alpha \neq \beta,\\ f(\alpha)=1}} \dot{A}_\alpha \cap \bigcap_{\substack{\alpha \neq \beta,\\ f(\alpha)=0}} (\omega \setminus \dot{A}_\alpha) \cap \dot{B}^{f(\beta)}$$

where $\dot{B}^1 = \dot{B}$ and $\dot{B}^0 = \omega \setminus \dot{B}$. By assumption on $\dot{B}$, we may find $(q', \bar{q}) \leq (p', \bar{p})$ with $F^{q'} \subseteq X$ and $k' \geq n^p$ such that

$$(*_6) \qquad (q', \bar{q}) \Vdash \forall\, i \geq k' \ \forall\, f \in 2^{\leq F^{p'} \cup \{\beta\}} \colon \left| \frac{|i \cap \dot{C}_f|}{i} - \frac{1}{2^{|\mathrm{dom}(f)|}} \right| < \frac{\varepsilon^p(f)}{16}.$$

[9] For $\kappa > \omega_1$, it is easy to see that $\mathbb{C}_\kappa$ still completely embeds into $\mathbb{P}_\kappa$ (this also follows from [BJZ97, Theorem 4.8], because $\mathbb{P}_\kappa$ is semi-Cohen), but not forcing-equivalent to $\mathbb{P}_\kappa$ any more.

We may assume $n^{q'} \geq k'$.

Now apply Lemma 4.3 with $p$, $E = X$, $m = k'$, $\varepsilon = \varepsilon^p$ and $q'$ to obtain $q$ such that $q \leq p$, $q \leq q'$, $F^q = F^{q'} \cup F^p$, $\varepsilon^q(f) = \varepsilon^p(f)$ for all $f \in 2^{\leq F^p} \setminus 2^{\leq F^{p'}}$, and $\varepsilon^q(f) \geq 16$ for all $f$ whose domain is not contained in either $F^p$ or $F^{q'}$. Let $q'' = q{\restriction}_{X \cup \{\beta\}}$. We may assume $q' = q{\restriction}_X = q''{\restriction}_X$.

Let $\ell \geq 8n^q$. We may find $(r', \bar{r}) \leq (q', \bar{q})$ with $F^{r'} \subseteq X$ such that $(r', \bar{r})$ decides $\dot{B} \cap \ell$. By Corollary 4.4, we may also assume

$$\frac{2^{2(|F^{r'}|+1)}}{n^{r'}} < \frac{\varepsilon^{r'}}{8} \tag{$*_7$}$$

and $n^{r'} \geq \ell$. Next, let $(s', \bar{s}) \leq (r', \bar{r})$ with $F^{s'} \subseteq X$ such that $(s', \bar{s})$ decides $\dot{B} \cap n^{r'}$. We now define a condition $r'' \in \mathbb{P}$ with $r'' \leq r'$ and $r'' \leq q''$ as follows:

- $F^{r''} = F^{r'} \cup \{\beta\} = F^{r'} \cup F^{q''}$,
- $n^{r''} = n^{r'}$,
- $a^{r''}_\alpha = a^{r'}_\alpha$ for $\alpha \in F^{r'}$, $a^{r''}_\beta \cap n^q = a^q_\beta$, and, for $n^q \leq i < n^{r'}$, $i \in a^{r''}_\beta$ iff $(s', \bar{s}) \Vdash i \in \dot{B}$, and
- $\varepsilon^{r''}{\restriction}_{2^{\leq F^{r'}}} = \varepsilon^{r'}$, $\varepsilon^{r''}(f) = \varepsilon^{q''}(f)$ for $f \in 2^{\leq F^{q''}}$ with $\beta \in \mathrm{dom}(f)$, and $\varepsilon^{r''}(f) \geq 16$ for all remaining $f$.

We need to check that $r''$ is indeed a condition and $r'' \leq q''$. ($r'' \leq r'$ then follows trivially.)

We first check (D5) for $r'' \leq q''$. Fix $i$ with $n^q \leq i \leq n^{r'}$. Also let $f \in 2^{\leq F^{p'} \cup \{\beta\}}$ with $\beta \in \mathrm{dom}(f)$. (There is nothing to show for other $f$, because they either belong to $2^{\leq F^{r'}}$ or they satisfy $\varepsilon^{r''}(f) \geq 16$.)

We will show only

$$\frac{|i \cap b^{r''}_f|}{i} < \frac{1}{2^{|\mathrm{dom}(f)|}} + \varepsilon^p(f),$$

since the other inequality is analogous. By assumption on $(q', \bar{q})$ and $(s', \bar{s})$, we know

$$(s', \bar{s}) \Vdash \left|n^q \cap \dot{C}_f\right| > n^q \cdot \left(\frac{1}{2^{|\mathrm{dom}(f)|}} - \frac{\varepsilon^p(f)}{16}\right)$$

and

$$(s', \bar{s}) \Vdash \left|i \cap \dot{C}_f\right| < i \cdot \left(\frac{1}{2^{|\mathrm{dom}(f)|}} + \frac{\varepsilon^p(f)}{16}\right).$$

Therefore

$$(s', \bar{s}) \Vdash \left|[n^q, i) \cap \dot{C}_f\right| < \frac{i - n^q}{2^{|\mathrm{dom}(f)|}} + \frac{n^q \cdot \varepsilon^p(f)}{16} + \frac{i \cdot \varepsilon^p(f)}{16}.$$

By the definition of $a^{r''}_\beta$, we now see that

$$\left|[n^q, i) \cap b^{r''}_f\right| < \frac{i - n^q}{2^{|\mathrm{dom}(f)|}} + \frac{n^q \cdot \varepsilon^p(f)}{16} + \frac{i \cdot \varepsilon^p(f)}{16}.$$

On the other hand, by (C5) for $q$ and $f$,

$$\left|n^q \cap b_f^{r''}\right| = \left|b_f^q\right| < n^q \cdot \left(\frac{1}{2^{|\mathrm{dom}(f)|}} + \frac{\varepsilon^p(f)}{8}\right).$$

Hence

$$\frac{\left|i \cap b_f^{r''}\right|}{i} < \frac{1}{2^{|\mathrm{dom}(f)|}} + \frac{n^q}{i} \cdot \frac{3 \cdot \varepsilon^p(f)}{16} + \frac{\varepsilon^p(f)}{16} < \frac{1}{2^{|\mathrm{dom}(f)|}} + \varepsilon^p(f),$$

as required for (D5). Furthermore, using $n^{r'} \geq 8n^q$, the previous formula with $i = n^{r'}$ gives

$$\frac{\left|n^{r'} \cap b_f^{r''}\right|}{n^{r'}} < \frac{1}{2^{|\mathrm{dom}(f)|}} + \frac{\varepsilon^p(f)}{8}$$

as required for (C5). On the other hand, since $|F^{r''}| = |F^{r'}| + 1$, (C6) is an immediate consequence of Eq. ($*_7$).

Finally, apply Lemma 4.3 with $p = q$, $E = X \cup \{\beta\}$, $p' = q''$, $m = \ell$, $\varepsilon = \varepsilon^q$ and $q' = r''$ to obtain $r$ with $r \leq q$, $r \leq r''$. In particular, we have $r \leq p$, and since $r \leq r''$, $(r, \bar{r})$ forces that $[n^q, \ell) \cap \dot{B} = [n^q, \ell) \cap \dot{A}_\beta$. Now note that

$$(r', \bar{r}) \Vdash \left|n^q \cap \dot{B}\right| < n^q \cdot \left(\frac{1}{2} + \frac{\varepsilon^p(f_0)}{16}\right)$$

and

$$(r', \bar{r}) \Vdash \left|\ell \cap \dot{B}\right| > \ell \cdot \left(\frac{1}{2} - \frac{\varepsilon^p(f_0)}{16}\right).$$

Therefore

$$(r, \bar{r}) \Vdash \left|[n^q, \ell) \cap \dot{B}\right| = \left|[n^q, \ell) \cap \dot{B} \cap \dot{A}_\beta\right| > \frac{\ell - n^q}{2} - \frac{\ell \cdot \varepsilon^p(f_0)}{8}$$

and hence, using $\ell \geq 8n^q$ and $\varepsilon^p(f_0) < 1/2$,

$$(r, \bar{r}) \Vdash \frac{\left|\ell \cap \dot{B} \cap \dot{A}_\beta\right|}{\ell} > \frac{\ell - n^q}{2\ell} - \frac{\varepsilon^p(f_0)}{8} > \frac{7}{16} - \frac{1}{16} = \frac{3}{8}$$

as required. □

Thus we obtain:

**Theorem 4.13.** *Let $\kappa$ be an uncountable cardinal. There is a generic extension with a maximal $1/2$-independent family of size $\kappa$.*

Using a finite support product of forcings $\mathbb{P}_\kappa$ together with an argument due to Blass ([Bla93, Theorem 9]), we see:

**Theorem 4.14.** *Let $V$ be a model of* ZFC *and* GCH. *In $V$, let $C$ be a closed set of uncountable cardinals with $\aleph_1 \in C$, $\kappa \in C$ for $\aleph_1 \leq \kappa \leq |C|$ and $\lambda^+ \in C$ for $\lambda \in C$ with $\mathrm{cf}(\lambda) = \omega$.*

*Then there is a ccc poset $\mathbb{Q}$ forcing $\mathfrak{c} = \max(C)$ and, in the generic extension, there is a maximal $1/2$-independent family of size $\kappa$ if and only if $\kappa \in C$.*

*Proof.* We force with the finite support product of the $\mathbb{P}_\kappa$ for $\kappa \in C$. Then by Lemma 4.12, there is a maximal $1/2$-independent family of size $\kappa$ for each $\kappa \in C$. The argument that there is no maximal $1/2$-independent family of size $\kappa$ for each $\kappa \notin C$ is exactly like the corresponding argument in [Bla93, Theorem 9] (see [BSZ00, Theorem 3.2] for a similar argument). □

Embedding the partial order $\mathbb{P}_\lambda$ (for $\lambda$ of countable cofinality) into the template framework as in [Bre03], we see:

**Theorem 4.15.** *Assume* CH *and let $\lambda$ be a singular cardinal of countable cofinality. Then there is a forcing extension satisfying $\mathfrak{i}_{1/2} = \lambda$. In particular, $\mathfrak{i}_{1/2} = \aleph_\omega$ is consistent.*

Roughly speaking, this can be proved by replacing Hechler forcing by localisation forcing (as in the proof of Lemma 3.9) and Hechler's poset for adding a mad family of size $\lambda$ by the poset $\mathbb{P}_\lambda$ in the framework of [Bre03]. Since many of the details are as in the latter article, we refrain from repeating the fifteen-page-long proof and only point out the differences. We ask the reader to have [Bre03] (and, as for the proof of Lemma 3.9, the recent [Bre21]) at hand.

Assume the linear order $\langle L, \leq_L \rangle$ is the disjoint union of the sets $L^{\text{prod}}$ and $L^{\text{iter}}$.[10] The intention is that $L^{\text{prod}}$ denotes the coordinates for generically adding a maximal $1/2$-independent family $\mathcal{A}$ with the forcing of Definition 4.1 and $L^{\text{iter}}$ is the set of coordinates for the $\mathbb{LOC}$-iterands exactly as in the proof of Lemma 3.9. The $\mathbb{LOC}$-generic with index $x \in L^{\text{iter}}$ will be generic over members of $\mathcal{A}$ with index $y < x$, $y \in L^{\text{prod}}$, but not over the others. To this end, we call a set $A \subseteq L$ *closed* iff for all $x \in A$ and $y < x$ with $y \in L^{\text{prod}}$, $y$ also belongs to $A$. For arbitrary $A \subseteq L$, its *closure* is defined by $\text{cl}(A) := A \cup \bigcup_{x \in A} L_x^{\text{prod}}$. We adapt the definition of "template" (Definition 3.10) to this context as follows:

**Definition 4.16** (see [Bre03, pp. 2635–2636])**.** A *template* is a pair $(L, \bar{\mathcal{I}})$ with $L = L^{\text{prod}} \cup L^{\text{iter}}$ and $\bar{\mathcal{I}} = \{\mathcal{I}_x \mid x \in L^{\text{iter}}\})$ such that $\langle L, \leq_L \rangle$ is a linear order, $\mathcal{I}_x \subseteq \mathcal{P}(L_x)$ consists of closed sets for $x \in L^{\text{iter}}$, and

(1) $\mathcal{I}_x$ contains $L_x^{\text{prod}}$ as well as $\text{cl}(\{y\})$ for all $y \in L_x^{\text{iter}}$ and is closed under unions and intersections,
(2) $\mathcal{I}_x \subseteq \mathcal{I}_y$ for $x <_L y$, and
(3) $\mathcal{I}{\restriction}_{L^{\text{iter}}} := \bigcup_{x \in L^{\text{iter}}} \mathcal{I}_x{\restriction}_{L^{\text{iter}}} \cup \{L^{\text{iter}}\}$ is well-founded with respect to inclusion, as witnessed by the *depth function* $\text{dp}_{\mathcal{I}{\restriction}_{L^{\text{iter}}}} : \mathcal{I}{\restriction}_{L^{\text{iter}}} \to \text{Ord}$.

If $A \subseteq L$ and $x \in L^{\text{iter}}$, we define $\mathcal{I}_x{\restriction}_A = \{B \cap A \mid B \in \mathcal{I}_x\}$ and let $\bar{\mathcal{I}}{\restriction}_A = \{\mathcal{I}_x{\restriction}_A \mid x \in A \cap L^{\text{iter}}\}$.

The depth function can easily be extended to all of $\mathcal{I} := \bigcup_{x \in L^{\text{iter}}} \mathcal{I}_x \cup \{L\}$ by $\text{dp}_{\mathcal{I}}(A) = \text{dp}_{\mathcal{I}{\restriction}_{L^{\text{iter}}}}(A \cap L^{\text{iter}})$ for $A \in \mathcal{I}$. This means a set has depth 0 iff it is a subset of $L^{\text{prod}}$. We can now recursively define the iteration, exactly as in [Bre21, Definition and Theorem 23], with the added twist that the basic step is not the trivial forcing but the forcing adding the maximal $1/2$-independent family. For closed $A \subseteq L$ define $\mathbb{P}{\restriction}_A$, basically as in [Bre03, Definition, pp. 2636–2637]:

[10] $L^{\text{prod}}$ and $L^{\text{iter}}$ are called $L_{\text{mad}}$ and $L_{\text{Hech}}$ in [Bre03], respectively.

- $\underline{\mathrm{dp}(A) = 0}$. This means $A \subseteq L^{\mathrm{prod}}$. Then $\mathbb{P}{\restriction}_A = \mathbb{P}_A$ in the sense of Definition 4.9 with the only change that $A$ is an arbitrary set and not necessarily a set of ordinals.
- $\underline{\mathrm{dp}(A) > 0}$. $\mathbb{P}{\restriction}_A$ consists of all finite partial functions $p$ with domain contained in $A$ such that $p{\restriction}_{L^{\mathrm{prod}}} \in \mathbb{P}_{A \cap L^{\mathrm{prod}}}$ and, letting $x = \max(\mathrm{dom}(p) \cap L^{\mathrm{iter}})$, there is $B \in \mathcal{I}_x{\restriction}_A$ such that $p{\restriction}_{A \cap L_x} \in \mathbb{P}{\restriction}_B$ and $p(x)$ is a $\mathbb{P}{\restriction}_B$-name for a condition in $\mathbb{L}\dot{\mathbb{O}}\mathbb{C}$.
  The ordering on $\mathbb{P}{\restriction}_A$ is given by $q \leq_{\mathbb{P}{\restriction}_A} p$ iff $\mathrm{dom}(q) \supseteq \mathrm{dom}(p)$, $q{\restriction}_{L^{\mathrm{prod}}} \leq_{\mathbb{P}_{L^{\mathrm{prod}}}} p{\restriction}_{L^{\mathrm{prod}}}$, and, letting $x = \max(\mathrm{dom}(q) \cap L^{\mathrm{iter}})$, there is $B \in \mathcal{I}_x{\restriction}_A$ such that $p{\restriction}_{A \cap L_x}, q{\restriction}_{A \cap L_x} \in \mathbb{P}{\restriction}_B$ and
  - either $x \notin \mathrm{dom}(p)$ and $q{\restriction}_{A \cap L_x} \leq_{\mathbb{P}{\restriction}_B} p{\restriction}_{A \cap L_x}$,
  - or $x \in \mathrm{dom}(p)$, $q{\restriction}_{A \cap L_x} \leq_{\mathbb{P}{\restriction}_B} p{\restriction}_{A \cap L_x}$, and $p(x)$ and $q(x)$ are $\mathbb{P}{\restriction}_B$-names for conditions in $\mathbb{L}\dot{\mathbb{O}}\mathbb{C}$ such that $q{\restriction}_{A \cap L_x} \Vdash_{\mathbb{P}{\restriction}_B} q(x) \leq_{\mathbb{L}\dot{\mathbb{O}}\mathbb{C}} p(x)$.

As remarked in both [Bre03] and [Bre21], to see that this definition works, one actually has to prove a number of facts along the iteration including complete embeddability. Most of this does not depend on the particular forcing notions iterated, but only on the template structure, and thus the arguments of [Bre03] carry over; however, there is one place in the proof of complete embeddability ([Bre03, Main Lemma 1.1]) where the iterands matter, and we therefore state the lemma and point out the changes.[11]

**Lemma 4.17.** *Let $B \in \mathcal{I}$ and $A \subset B$ be closed. Then $\mathbb{P}{\restriction}_B$ is a partial order, $\mathbb{P}{\restriction}_A \subseteq \mathbb{P}{\restriction}_B$ and even $\mathbb{P}{\restriction}_A \lessdot \mathbb{P}{\restriction}_B$. More explicitly, any $p \in \mathbb{P}{\restriction}_B$ has a* canonical reduction *$p_0 = p_0(p, A, B) \in \mathbb{P}{\restriction}_A$ such that*

*(1)* $\mathrm{dom}(p_0) = \mathrm{dom}(p) \cap A$,
*(2)* $\sigma_x^{p_0} = \sigma_x^p$ *for all* $x \in \mathrm{dom}(p_0) \cap L^{\mathrm{iter}}$ *and* $p_0(x) = p(x)$ *for all* $x \in \mathrm{dom}(p_0) \cap L^{\mathrm{prod}}$,

*and such that, whenever $D \in \mathcal{I}$, $B, C \subseteq D$, $C$ closed, $C \cap B = A$, and $q_0 \in \mathbb{P}{\restriction}_C$ extends $p_0$, then there is $q \in \mathbb{P}{\restriction}_D$ extending both $q_0$ and $p$.*

Recall here that for $x \in \mathrm{dom}(p) \cap L^{\mathrm{iter}}$, $p(x) = (\sigma_x^p, \dot{\varphi}_x^p)$ is a $\mathbb{P}{\restriction}_{\bar{B}}$-name for a condition in $\mathbb{L}\dot{\mathbb{O}}\mathbb{C}$ for some $\bar{B} \in \mathcal{I}_x{\restriction}_B$, and we may assume that $p{\restriction}_{\bar{B}}$ decides the first coordinate (because it is finite), that is, $\sigma_x^p$ is not a name.

*Proof.* We conduct recursion induction on $\mathrm{dp}(B) = \alpha$. Work with [Bre03, Main Lemma 1.1] at hand. The parts referring to the template (and this is most of the proof) are exactly as in that proof, so we will not repeat them. However, since $\mathbb{L}\mathbb{O}\mathbb{C}$ is different from $\mathbb{D}$, the construction of the *canonical projection* $\dot{\varphi}_x^{p_0}$ of $\dot{\varphi}_x^p$ (bottom half of p. 2638 in [Bre03]) has to be changed, as follows:

We are given $p \in \mathbb{P}{\restriction}_B$ and $x = \max(\mathrm{dom}(p) \cap L^{\mathrm{iter}})$. There is $\bar{B} \in \mathcal{I}_x{\restriction}_B$ such that $\bar{p} = p{\restriction}_{B \cap L_x} \in \mathbb{P}{\restriction}_{\bar{B}}$ and $\dot{\varphi}_x^p$ is a $\mathbb{P}{\restriction}_{\bar{B}}$-name. Letting $\bar{A} = A \cap \bar{B}$, we see $\bar{A} \in \mathcal{I}_x{\restriction}_A$. By induction hypothesis, $\mathbb{P}{\restriction}_{\bar{A}} \lessdot \mathbb{P}{\restriction}_{\bar{B}}$ and $\bar{p}$ has a reduction $\bar{p}_0 = p_0(\bar{p}, \bar{A}, \bar{B}) \in \mathbb{P}{\restriction}_{\bar{A}}$. We may assume $x \in A$, for otherwise, there is nothing to show (Case 2 in [Bre03]).

[11] An alternative approach would be to redo [Bre21, Definition and Theorem 23] in a framework with $L^{\mathrm{prod}}$. This is more general, but also involves more work. Therefore we stick to [Bre03].

Work with the cBa's $\mathbb{B}_{\bar A} = \mathrm{r.o.}(\mathbb{P}\restriction_{\bar A})$ and $\mathbb{B}_{\bar B} = \mathrm{r.o.}(\mathbb{P}\restriction_{\bar B})$. Let $n := |\sigma^p_x|$. Thus $\bar p \Vdash |\dot\varphi^p_x(i)| \le n$ for all $i$. Consider a partial function $\tau$ such that $\mathrm{dom}(\tau) = [n, m)$ for some $m = m_\tau \ge n$ and $\tau(i) \in [\omega]^{\le n}$ for all $i \in \mathrm{dom}(\tau)$, and let $b_\tau = [\![\forall i \in \mathrm{dom}(\tau)\colon \dot\varphi^p_x(i) = \tau(i)]\!] \cap \bar p$. Notice that for fixed $m > n$, the $b_\tau$ with $m_\tau = m$ form a maximal antichain below $\bar p$. Let $a^*_\tau$ be the product of $\bar p_0$ and the projection of $b_\tau$ to $\mathbb{B}_{\bar A}$. Then $\sum\{a^*_\tau \mid m_\tau = m\} = \bar p_0$ for $m > n$. Define $a_\tau$ by recursion on $m_\tau \ge n$ as follows: $a_\tau = \bar p_0$ for $m_\tau = n$ (note that $\tau$ is the empty function in this case). Let $\{u_j \mid j \in \omega\}$ list $[\omega]^{\le n}$. For $\tau$ with $m = m_\tau > n$ let $a_\tau = a_{\tau\restriction_{m-1}} \cdot (a^*_\tau \setminus \sum_{j<k} a_{\tau\restriction_{m-1}{}^\frown\langle u_j\rangle})$ where $k$ is unique such that $\tau(m-1) = u_k$. It is easy to see that for fixed $m > n$, the $a_\tau$ with $m_\tau = m$ form a maximal antichain below $\bar p_0$. Therefore they canonically define a $\mathbb{P}\restriction_{\bar A}$-name $\dot\varphi^{p_0}_x$ such that $\bar p_0 \Vdash_{\mathbb{P}\restriction_{\bar A}} |\dot\varphi^{p_0}_x(i)| \le n$ for all $i$.[12]

The main property of this name is that for all $\sigma \in ([\omega]^{<\omega})^{<\omega}$ with $|\sigma(i)| \le i$ for $i \in |\sigma|$ and $\sigma^p_x \subseteq \sigma$, $a'_\sigma = \sum\{a_\tau \mid m_\tau = |\sigma| \text{ and } \tau(i) \subseteq \sigma(i) \text{ for all } i \text{ with } n \le i < |\sigma|\}$ is a reduction of $b'_\sigma = \sum\{b_\tau \mid m_\tau = |\sigma| \text{ and } \tau(i) \subseteq \sigma(i) \text{ for all } i \text{ with } n \le i < |\sigma|\}$.

The rest of the argument can now be completed as in [Bre03]. When defining the stronger condition $q$ on coordinates from $L^{\mathrm{prod}}$ (see the two cases on pp. 2639–2640 in [Bre03]), use Lemma 4.3. □

The whole forcing $\mathbb{P}\restriction_L$ is ccc [Bre03, Lemma 1.2], and if $\mu$ is regular uncountable, $\mu \subseteq L^{\mathrm{iter}}$ is cofinal in $L$, and $L_\alpha \in \mathcal{I}_\alpha$ for all $\alpha \in \mu$, then $\mathbb{P}\restriction_L$ forces $\mathrm{add}(\mathcal{N}) = \mathrm{cof}(\mathcal{N}) = \mu$ (this is like [Bre03, Proposition 1.6] and exactly as in the third paragraph of the proof of Lemma 3.9).

**Lemma 4.18.** *Assume $L$ has uncountable cofinality and $L^{\mathrm{prod}}$ is cofinal in $L$. Then $\mathbb{P}\restriction_L$ adds a maximal $1/2$-independent family.*

*Proof.* This is analogous to [Bre03, Proposition 1.7]. However, since $\mathbb{P}_{L^{\mathrm{prod}}}$ is a much more complicated forcing than the forcing adding a mad family, the proof is more complex and the combinatorial and computational details are much more like the proof of Lemma 4.12 above. We provide an outline, explaining in detail how to adapt the latter proof to the present context.

Let $\mathcal{A} = \{A_x \mid x \in L^{\mathrm{prod}}\}$ be the generic $1/2$-independent family added by $\mathbb{P}\restriction_{L^{\mathrm{prod}}}$ (which completely embeds into $\mathbb{P}\restriction_L$). See between Corollary 4.10 and Corollary 4.11 for the definition. We need to check maximality, so let $\dot B$ be a $\mathbb{P}\restriction_L$-name for an infinite and coinfinite subset of $\omega$. By [Bre03, Lemma 1.4] (cf. [Bre21, Lemma 25]), there is a countable set $C \subseteq L$ such that $\dot B$ is a $\mathbb{P}\restriction_{\mathrm{cl}(C)}$-name. Since $L$ has uncountable cofinality and $L^{\mathrm{prod}}$ is cofinal in $L$, there is $x \in L^{\mathrm{prod}}$ such that $\mathrm{cl}(C) \subseteq L_x$. Therefore $\dot B$ is a $\mathbb{P}\restriction_{L_x}$-name.

**Claim.** *Assume $\tilde p_0 \in \mathbb{P}\restriction_L$ forces that $\dot B$ is $1/2$-independent from all $\dot A_y$ for $y \in L^{\mathrm{prod}}_x$. Then $\tilde p_0$ forces that for all $k$, there is an $\ell > k$ such that*

$$\frac{|\ell \cap \dot B \cap \dot A_x|}{\ell} > \frac{3}{8}.$$

[12] Notice that the values of $\dot\varphi^p_x$ (and $\dot\varphi^{p_0}_x$) at $i < n$ are irrelevant for the definition of the partial order, and that is why we omitted them here; one may actually assume that the domain of these functions is $[n, \omega)$.

Fix $\tilde{p} \leq \tilde{p}_0$ in $\mathbb{P}\upharpoonright_L$ and $k$. We need to find $\ell > k$ and $\tilde{r} \leq \tilde{p}$ in $\mathbb{P}\upharpoonright_L$ forcing the required statement. Let $\check{p} = \tilde{p}\upharpoonright_{L_x \cup L^{\mathrm{prod}}}$. Now redo the proof of Lemma 4.12 with $\check{p}$ instead of $(p, \bar{p})$ (and similarly $\check{q}$, $\check{r}$, ... instead of $(q, \bar{q})$, $(r, \bar{r})$, ...). Furthermore, let $p = \tilde{p}\upharpoonright_{L^{\mathrm{prod}}}$, $p' = \tilde{p}\upharpoonright_{L_x^{\mathrm{prod}}}$ and $\hat{p} = \tilde{p}\upharpoonright_{L_x}$ (and similarly for $q$, $r$, ...). $\hat{p}$ plays the role of $(p', \bar{p})$ and $x$ and $L_x^{\mathrm{prod}}$ play the roles of $\beta$ and $X$, respectively. Choices for other items are exactly like in Lemma 4.12.

As before, first find $\hat{q} \leq \hat{p}$ in $\mathbb{P}\upharpoonright_{L_x}$ and $k'$ such that $\hat{q}$ forces the statement in Eq. ($*_6$). From Lemma 4.3 we obtain again $q$ extending both $q' = \hat{q}\upharpoonright_{L^{\mathrm{prod}}}$ and $p$ in $\mathbb{P}\upharpoonright_{L^{\mathrm{prod}}}$. This gives us $\check{q}$ in $\mathbb{P}\upharpoonright_{L_x \cup L^{\mathrm{prod}}}$ such that $\check{q}\upharpoonright_{L^{\mathrm{prod}}} = q$ and $\check{q}\upharpoonright_{L^{\mathrm{iter}}} = \hat{q}\upharpoonright_{L^{\mathrm{iter}}}$. Also let $q'' = q\upharpoonright_{L_x \cup \{x\}}$.

Next, fix $\ell$ and find $\hat{r} \leq \hat{q}$ in $\mathbb{P}\upharpoonright_{L_x}$ deciding $\dot{B} \cap \ell$, and then $\hat{s} \leq \hat{r}$ deciding $\dot{B} \cap n^{r'}$ (where $r' = \hat{r}\upharpoonright_{L^{\mathrm{prod}}}$). Define $r'' \leq r', q''$ in $\mathbb{P}\upharpoonright_{L_x^{\mathrm{prod}} \cup \{x\}}$ as before. The proof that this is indeed a condition extending $q''$ carries over verbatim. Thus we can again apply Lemma 4.3 to obtain $r$ extending both $r''$ and $q$ in $\mathbb{P}\upharpoonright_{L^{\mathrm{prod}}}$. We then see that $\check{r} \in \mathbb{P}\upharpoonright_{L_x \cup L^{\mathrm{prod}}}$, defined by $\check{r}\upharpoonright_{L^{\mathrm{prod}}} := r$ and $\check{r}\upharpoonright_{L^{\mathrm{iter}}} := \hat{r}\upharpoonright_{L^{\mathrm{iter}}}$, forces the required statement.

In the last step (note that this is different from the proof of Lemma 4.12), we define $\tilde{r} \in \mathbb{P}\upharpoonright_L$ by $\tilde{r}\upharpoonright_{L_x \cup L^{\mathrm{prod}}} := \check{r}$ and $\tilde{r}(y) := \tilde{p}(y)$ for $y \in \mathrm{dom}(p) \cap (L^{\mathrm{iter}} \setminus L_x)$. Clearly $\tilde{r} \leq \tilde{p}, \check{r}$ is as required. □

We next recall the definition of the template we are using here [Bre03, Definition, p. 2643]. Assume $\lambda_0 \geq \aleph_2$ is regular and $\lambda > \lambda_0$ is a singular cardinal of countable cofinality, say $\lambda = \bigcup_n \lambda_n$ with $\lambda_n$ regular, equal to $\lambda_n^{\aleph_0}$ and strictly increasing. Also suppose $\kappa^{\aleph_0} < \lambda_n$ for $\kappa < \lambda_n$. For each $n$, choose a partition $\lambda_n^* = \bigcup_{\alpha<\omega_1} S_n^\alpha$ such that each $S_n^\alpha$ is coinitial in $\lambda_n^*$. Also assume $S_n^\alpha \cap \lambda_m^* = S_m^\alpha$ for $m < n$.

**Definition 4.19** (cf. Definition 3.11)**.** Elements of $L$ are non-empty finite sequences $x$ such that

- $x(0) \in \lambda_0$,
- $x(n) \in \lambda_n^* \cup \lambda_n$ for $0 < n < |x| - 1$, and
- in case $|x| \geq 2$,
  - if $x(|x|-2)$ is positive, then $x(|x|-1) \in \lambda_{|x|-1}^* \cup \lambda$,
  - and if $x(|x|-2)$ is negative, then $x(|x|-1) \in \lambda^* \cup \lambda_{|x|-1}$.

Let $x \in L^{\mathrm{iter}}$ if $|x| = 1$ or $x(|x|-1) \in \lambda_{|x|-1}^* \cup \lambda_{|x|-1}$; otherwise $x \in L^{\mathrm{prod}}$. The order on $L$ is naturally given by $x < y$ if

- either $x \supsetneq y$ and $y(|x|)$ is positive,
- or $y \supsetneq x$ and $x(|y|)$ is negative,
- or, letting $n := \min\{m \mid x(m) \neq y(m)\}$, $x(n) < y(n)$ in the natural ordering of $\lambda^* \cup \lambda$.

Clearly, this is a linear ordering. Identify sequences of length one with their ranges so that $\lambda_0$ is cofinal in $L$. "Relevant" members of $L^{\mathrm{iter}}$ are defined exactly as before, between Definition 3.11 and Definition 3.12 (see also [Bre03, p. 2643]).

**Definition 4.20** (cf. Definition 3.12)**.** For $x \in L^{\mathrm{iter}}$, let $\mathcal{I}_x$ consist of finite unions of

- $L_\alpha$ where $\alpha \leq x$ and $\alpha \in \lambda_0$,
- $\mathrm{cl}(J_y)$ where $y \leq x$ is relevant,
- $\mathrm{cl}(\{y\})$ where $y \in L_x^{\mathrm{iter}}$, and
- $L_x^{\mathrm{prod}}$.

Then $(L, \bar{\mathcal{I}})$ with $\mathcal{I} = \{\mathcal{I}_x \mid x \in L^{\mathrm{iter}}\}$ is indeed a template [Bre03, Lemma 2.1].

*Proof of Theorem 4.15.* Take the template $(L, \bar{\mathcal{I}})$ introduced above and let $\mathbb{P} = \mathbb{P}{\restriction}_L$. By the properties of $L$, we see that $\mathbb{P}$ forces $\mathrm{add}(\mathcal{N}) = \mathrm{cof}(\mathcal{N}) = \lambda_0$ (see the comment immediately preceding Lemma 4.18) and that it adds a maximal $1/2$-independent family of size $\lambda$ by Lemma 4.18. Thus $\mathfrak{r}_{1/2} = \lambda_0 \leq \mathfrak{i}_{1/2} \leq \lambda$ (cf. [Bre03, Corollary 2.2]). Therefore it suffices to show that there are no maximal $1/2$-independent families of size $\kappa$ with $\lambda_0 \leq \kappa < \lambda$. This is done by an isomorphism-of-names argument; since this argument does not depend on the individual forcings, but only on the structure of the template, it works exactly as in the corresponding proof in [Bre03, Section 3]. Only the very end of this proof needs to be changed, in a way similar to how the proof of [Bre21, Theorem 30] was changed in Lemma 3.9. We provide the details of this last step and request the reader again to have [Bre03] at hand.

In the last paragraph of that proof of the latter (p. 2648), we know by construction that whenever $F \subset \kappa$ is finite, we can find $\alpha < \omega_1$ such that $B^\kappa \cup \bigcup_{\beta \in F} B^\beta$ and $B^\alpha \cup \bigcup_{\beta \in F} B^\beta$ are weakly isomorphic via the mapping fixing nodes of $\bigcup_{\beta \in F} B^\beta$ and sending the $x_s^\kappa$ to the corresponding $x_s^\alpha$, and such that this mapping identifies cofinal subsets of the traces of $\mathcal{I}$ on the two sets (in fact, this is true for all but countably many $\alpha$). This weak isomorphism canonically extends to a weak isomorphism of $C^\kappa \cup \bigcup_{\beta \in F} C^\beta$ and $C^\alpha \cup \bigcup_{\beta \in F} C^\beta$, which in turn means that $\mathbb{P}{\restriction}_{C^\kappa \cup \bigcup_{\beta \in F} C^\beta}$ and $\mathbb{P}{\restriction}_{C^\alpha \cup \bigcup_{\beta \in F} C^\beta}$ are isomorphic [Bre03, Lemma 3.2] by a mapping sending $\dot{A}^\kappa$ to $\dot{A}^\alpha$. Since $\{\dot{A}^\beta \mid \beta \in \{\alpha\} \cup F\}$ is forced to be $1/2$-independent (by $\mathbb{P}{\restriction}_{C^\alpha \cup \bigcup_{\beta \in F} C^\beta}$), $\{\dot{A}^\beta \mid \beta \in \{\kappa\} \cup F\}$ is forced to be $1/2$-independent (by the isomorphic $\mathbb{P}{\restriction}_{C^\kappa \cup \bigcup_{\beta \in F} C^\beta}$). This completes the proof of the non-maximality of a family $\dot{\mathcal{A}}$ of size $\kappa$. □

For a similar argument, cf. [FT15]. Note that since $\mathrm{cov}(\mathcal{N})$ is a lower bound of $\mathfrak{i}_{1/2}$, it is clear (and much easier to prove) that $\mathfrak{i}_{1/2}$ can be a singular cardinal of uncountable cofinality (in the appropriate random model).

## 5. Open Questions

While we have shown that several of our newly defined cardinal characteristics are, in fact, new, there are still a number of open questions.

**Question A.** *We summarise the open questions related to Figure 1:*

*(Q1) Does* $\mathrm{Con}(\mathfrak{s}_{1/2} < \mathrm{non}(\mathcal{N}))$ *hold or is* $\mathfrak{s}_{1/2} = \mathrm{non}(\mathcal{N})$*?*

*(Q2) Does* $\mathrm{Con}(\mathfrak{d} < \mathfrak{s}_{1/2\pm\varepsilon} \leq \mathfrak{s}_{1/2})$ *hold or is* $\mathfrak{s}_{1/2} \leq \mathfrak{d}$*? (If it is the latter, we already know* $\mathrm{Con}(\mathfrak{s}_{1/2} < \mathfrak{d})$ *by* $\mathrm{Con}(\mathrm{non}(\mathcal{N}) < \mathfrak{d})$*.)*[13]

[13] Since the writing of this paper, the consistency of $\mathfrak{d} < \mathfrak{s}_{1/2\pm\varepsilon}$ as well as of the dual $\mathfrak{b} > \mathfrak{r}_{1/2\pm\varepsilon}$ has been proved independently by Farkas, Klausner and Lischka and by Valderrama.

(Q3) *Which of the following statements are true?*

$$\mathrm{Con}(\mathfrak{s} < \mathfrak{s}^{w}_{1/2}) \quad \textit{or} \quad \mathfrak{s} = \mathfrak{s}^{w}_{1/2}$$
$$\mathrm{Con}(\mathfrak{s}^{w}_{1/2} < \mathfrak{s}^{\infty}_{1/2}) \quad \textit{or} \quad \mathfrak{s}^{w}_{1/2} = \mathfrak{s}^{\infty}_{1/2}$$
$$\mathrm{Con}(\mathfrak{s}_{1/2\pm\varepsilon} < \mathfrak{s}_{1/2}) \quad \textit{or} \quad \mathfrak{s}_{1/2\pm\varepsilon} = \mathfrak{s}_{1/2}$$

(Q4) *Given $\varepsilon > \varepsilon'$ and an $\varepsilon$-almost bisecting family, can one (finitarily) modify it to get an $\varepsilon'$-almost bisecting family of equal size? (If yes, then $\mathfrak{s}_{1/2\pm\varepsilon}$ is independent of $\varepsilon$. If not, then $\inf_{\varepsilon\in(0,1/2)} \mathfrak{s}_{1/2\pm\varepsilon}$ and $\sup_{\varepsilon\in(0,1/2)} \mathfrak{s}_{1/2\pm\varepsilon}$ might be interesting characteristics, as well.)*

(Q5) *Can characteristics in the upper row of the diagram consistently be smaller than ones in the lower row? Specifically, which of the following statements are true?*

$$\mathrm{Con}(\mathfrak{s}_{1/2\pm\varepsilon} < \mathfrak{s}^{w}_{1/2}) \quad \textit{or} \quad \mathfrak{s}_{1/2\pm\varepsilon} \geq \mathfrak{s}^{w}_{1/2}$$
$$\mathrm{Con}(\mathfrak{s}_{1/2\pm\varepsilon} < \mathfrak{s}^{\infty}_{1/2}) \quad \textit{or} \quad \mathfrak{s}_{1/2\pm\varepsilon} \geq \mathfrak{s}^{\infty}_{1/2}$$
$$\mathrm{Con}(\mathfrak{s}_{1/2} < \mathfrak{s}^{\infty}_{1/2}) \quad \textit{or} \quad \mathfrak{s}_{1/2} \geq \mathfrak{s}^{\infty}_{1/2}$$

**Question B.** *We summarise the open questions related to Figure 2:*

(Q6) *Is it consistent that $\mathfrak{i}_* < 2^{\aleph_0}$?*
(Q7) *Which relations between $\mathfrak{i}_{1/2}$, $\mathfrak{i}_*$ and $\mathfrak{i}$ are true or consistent?*
(Q8) *Are there any smaller upper bounds for $\mathfrak{i}_{1/2}$ and $\mathfrak{i}_*$?*
(Q9) *Which relations between $\mathfrak{s}_{1/2}$ and $\mathfrak{s}_*$ are true or consistent?*
(Q10) *Which of the following statements are true?*

$$\mathrm{Con}(\mathrm{cov}(\mathcal{N}) < \mathfrak{r}_{1/2}) \quad \textit{or} \quad \mathrm{cov}(\mathcal{N}) = \mathfrak{r}_{1/2}$$
$$\mathrm{Con}(\mathfrak{r}_{1/2} < \mathfrak{r}_*) \quad \textit{or} \quad \mathfrak{r}_{1/2} = \mathfrak{r}_*$$
$$\mathrm{Con}(\mathfrak{s}_* < \mathrm{non}(\mathcal{N})) \quad \textit{or} \quad \mathfrak{s}_* = \mathrm{non}(\mathcal{N})$$

We suspect that (Q6) might be provable (via $\mathrm{Con}(\mathfrak{i}_* < \mathfrak{i})$) using the same idea as in Lemma 3.16. In an earlier draft of this paper, we had a somewhat convoluted creature forcing argument (with the help of some probabilistic sleight of hand) intended to prove $\mathrm{Con}(\mathfrak{s}_{1/2} < \mathrm{non}(\mathcal{N}))$, which unfortunately turned out to be incorrect. It seems plausible that such a creature forcing construction might be able to prove the intended result, after all; if an analogous probabilistic argument can be reproduced for $\mathfrak{s}_*$, a similar approach might also work to answer the third part of (Q10) and prove $\mathrm{Con}(\mathfrak{s}_* < \mathrm{non}(\mathcal{N}))$. Finally, since it is not too difficult to ensure that a creature forcing poset keeps $\mathrm{cov}(\mathcal{N})$ small (compare [FGKS17, Lemma 5.4.2] or [GK21, Lemma 7.8]), a clever creature forcing construction might be able to answer the first part of (Q10) and prove $\mathrm{Con}(\mathrm{cov}(\mathcal{N}) < \mathfrak{r}_{1/2})$.

Graduate School of System Informatics, Kobe University, 1-1 Rokkodai-cho, Nada-ku, Kobe, Hyogo Prefecture, 657-8501, Japan

*Email address*: brendle@kobe-u.ac.jp

Department of Mathematics, ETH Zurich, Rämistrasse 101, 8092 Zurich, Switzerland

*Email address*: lorenz.halbeisen@math.ethz.ch

Institute of Discrete Mathematics and Geometry, TU Wien, Wiedner Hauptstrasse 8–10/104, 1040 Wien, Austria

*Email address*: mail@l17r.eu

*URL*: https://l17r.eu

Department of Mathematics, ETH Zurich, Rämistrasse 101, 8092 Zurich, Switzerland

*Email address*: marc.lischka@gmail.com

Einstein Institute of Mathematics, Edmond J. Safra Campus, The Hebrew University of Jerusalem, Givat Ram, Jerusalem, 9190401, Israel and Department of Mathematics, Hill Center – Busch Campus, Rutgers, The State University of New Jersey, 110 Frelinghuysen Road, Piscataway, NJ 08854-8019, United States

*Email address*: shelah@math.huji.ac.il